\documentclass[12pt]{article}
\usepackage{epsfig,wrapfig,epic}
\usepackage{amsmath,latexsym}
\usepackage{amscd}
\usepackage{mathtools}
\usepackage{rotating}
\usepackage{soul}
\usepackage{fancybox}
\usepackage[normalem]{ulem}

\usepackage{hyperref}
\usepackage{amsfonts}
\usepackage{mathrsfs}
\usepackage{stmaryrd,bbm}
\usepackage{shadow}
\usepackage{lineno}
\usepackage{setspace}
\usepackage{color}
\newcommand{\blue}[1]{\textcolor{blue}{#1}}

\usepackage{dsfont}
\usepackage{graphics}
\usepackage{graphicx}
\usepackage{amssymb}

\title{On the Sensitivity of Singular and Ill-Conditioned Linear Systems}
\author{Zhonggang Zeng \thanks{Department of Mathematics,
Northeastern Illinois University, Chicago, Illinois 60625, USA.
~~email:~{\tt zzeng@neiu.edu}. 
~Research is supported in part by NSF under grant DMS-1620337.}}

\DeclareMathSymbol{\bdC}{\mathbin}{AMSb}{'103}
\DeclareMathAlphabet{\mathpzc}{OT1}{pzc}{m}{it}
\DeclareMathSymbol{\mP}{\mathbin}{AMSb}{'120}

\begin{document}
\newcommand{\vvv}{|\! |\! |}
\newcommand{\bdu}{\mathbf{u}}
\newcommand{\bdv}{\mathbf{v}}
\newcommand{\bdw}{\mathbf{w}}
\newcommand{\bdx}{\mathbf{x}}
\newcommand{\bdy}{\mathbf{y}}
\newcommand{\bdz}{\mathbf{z}}
\newcommand{\bda}{\mathbf{a}}
\newcommand{\bdb}{\mathbf{b}}
\newcommand{\bdc}{\mathbf{c}}
\newcommand{\bdd}{\mathbf{d}}
\newcommand{\bde}{\mathbf{e}}
\newcommand{\bdf}{\mathbf{f}}
\newcommand{\bdg}{\mathbf{g}}
\newcommand{\bdh}{\mathbf{h}}
\newcommand{\bdi}{\mathbf{i}}
\newcommand{\bdj}{\mathbf{j}}
\newcommand{\bdm}{\mathbf{m}}
\newcommand{\bdn}{\mathbf{n}}
\newcommand{\bdk}{\mathbf{k}}
\newcommand{\bdl}{\mathbf{l}}
\newcommand{\bdr}{\mathbf{r}}
\newcommand{\bdp}{\mathbf{p}}
\newcommand{\bdq}{\mathbf{q}}
\newcommand{\bds}{\mathbf{s}}
\newcommand{\bdo}{\mathbf{0}}
\newcommand{\bdF}{\mathbf{F}}
\newcommand{\C}{\mathbbm{C}}
\newcommand{\cK}{\mathpzc{Kernel}}
\newcommand{\cR}{\mathpzc{Range}}
\newcommand{\rank}[1]{\mathpzc{rank}\left(\,#1\,\right)}
\newcommand{\ranka}[2]{\mathpzc{rank}_{#1}\left(\,#2\,\right)}
\newcommand{\dist}[1]{\mathpzc{dist}\left(\,#1\,\right)}
\newcommand{\h}{{{\mbox{\tiny $\mathsf{H}$}}}}
\newcommand{\tr}{{{\mbox{\tiny $\top$}}}}
\newcommand{\sg}{\sigma}
\newcommand{\spn}{\mathpzc{span}}
\newcommand{\cG}{{\cal G}}
\newcommand{\cA}{{\cal A}}
\newcommand{\cV}{{\cal V}}
\newcommand{\Dl}{\Delta}
\newcommand{\bt}{\beta}
\newcommand{\cD}{{\cal D}}
\newcommand{\al}{\alpha}
\newcommand{\eps}{\varepsilon}
\newcommand{\la}{\lambda}
\newcommand{\cM}{{\cal M}}
\newcommand{\cU}{{\cal U}}
\newcommand{\dl}{\delta}
\newcommand{\Qed}{${~} $ \hfill \raisebox{-0.3ex}{\LARGE $\Box$}}

\newcommand{\sol}{\mathpzc{sol}}
\newtheorem{algrthm}{Algorithm}
\newtheorem{example}{Example}
\newtheorem{problem}{Problem}
\newtheorem{remark}{Remark}
\newtheorem{theorem}{Theorem}
\newtheorem{corollary}{Corollary}
\newtheorem{lemma}{Lemma}
\newtheorem{definition}{Definition}

\input amssym.def
\maketitle


\begin{abstract}
~Solving a singular linear system for an individual vector solution is an 
ill-posed problem with a condition number infinity.
~From an alternative perspective, however, the general solution of a 
singular system is of a bounded sensitivity as a unique element in an affine 
Grassmannian.
~If a singular linear system is given through empirical data that
are sufficiently accurate with a tight error bound, a properly formulated 
general numerical solution uniquely exists in the same affine 
Grassmannian, enjoys Lipschitz continuity and approximates the 
underlying exact solution with an accuracy in the same order as the 
data.
~Furthermore, any backward accurate numerical solution vector is an accurate
approximation to one of the solutions of the underlying singular system.
\end{abstract}

\section{Introduction}

Solving linear systems in the matrix-vector form 
~$A\,\bdx\,=\,\bdb$ ~is one of the most fundamental 
problems in scientific computing.
~In the literature of numerical analysis, 
linear systems are always assumed to be nonsingular with few exceptions.
~Numerical solutions of singular 
systems are almost never mentioned directly in textbooks. 
~A rare remark in Meyer's textbook \cite[page 218]{Meyer} accurately reflects the
state of knowledge:
~``If ~$A$ ~is singular, ... even a stable algorithm can result in a 
significant loss of information. ... 
[T]he small perturbation ~$E$ ~due to roundoff 
makes the possibility that ~$\rank{A+E}\,>\,\rank{A}$ ~very likely. 
~{\em The moral is to avoid floating point solutions of singular systems}''
(emphasis added).
~In applications such as deblurring images 
and discrete inverse problems, rank-deficient and 
highly ill-conditioned linear systems are approached using the 
Tikhonov regularization \cite{HansenBook,HansInv,HanNagOle,Neu98}.
~As Neumaier states \cite{Neu98}:
~``Though frequently needed in applications, the adequate handling of 
such ill-posed linear problems is hardly ever touched upon in numerical
analysis text books.''

Singular linear systems are unavoidable in scientific computing and often
need to be solved without knowing the exact matrices and vectors, as shown
in case studies in \S\ref{s:mod}. 
\,The obvious difficulty in solving a singular linear system from empirical
data is the condition number infinity so that the error is unbounded when 
solving for an individual vector solution.
~While this error analysis in itself is impeccable, the solution of a singular
system is more than an individual vector.
~The very notion of the numerical solution to a given system 
~$\tilde{A}\,\bdx\,=\,\tilde\bdb$ ~needs clarification when entries of
~$(\tilde{A},\tilde\bdb)$ ~serve as empirical data for an underlying singular
linear system ~$A\,\bdx\,=\bdb$.

This paper attempts to analyze the accuracy and sensitivity of solving 
singular linear systems from a different perspective: 
~The solution of a singular linear system is either an empty set or an 
affine subspace as a unique element in an affine Grassmannian rather than a 
vector.
~Using this point of view, the condition number becomes bounded.
~A properly formulated general numerical solution in a certain affine 
Grassmannian 
is of a sensitivity proportional to ~$\|A\|_2\,\|A^\dagger\|_2$, ~never 
infinity, with respect to either constrained or arbitrary perturbations
where ~$A^\dagger$ ~is the Moore-Penrose inverse of ~$A$.
~Such a numerical solution of a perturbed system 
~$\tilde{A}\,\bdx\,=\,\tilde\bdb$ ~within a viable error tolerance 
accurately solves the underlying singular system
~$A\,\bdx\,=\,\bdb$ ~and the ratio of solution accuracy to the data error is 
bounded by a factor of ~$\|A\|_2\,\|A^\dagger\|_2$, ~not
~$\|\tilde{A}\|_2\,\|\tilde{A}^{-1}\|_2$, 
~assuming the data error is small with an attainable tight bound.

We shall further demonstrate that the sensitivity of a singular 
linear system ~$A\,\bdx\,=\,\bdb$ ~is measured by
~$\|A\|_2\,\|A^\dagger\|_2$ ~rather than infinity from multiple 
perspectives, including homogeneous cases, under constrained perturbations
preserving the singularity and consistency, solving for the 
general numerical solutions in an affine Grassmannian, and solving for a 
single particular solution.
~Furthermore, every backward accurate numerical (vector) solution of a singular
consistent linear system accurately approximates a particular exact solution
regardless of the algorithm used.
~The ``error'' largely falls harmlessly in the kernel of ~$A$.
~This result extends what Peters and Wilkinson discovered in \cite{PetWil79} 
beyond inverse power iterations.
~While any numerical (single-vector) solution may be inaccurate to a 
linear system that is genuinely nonsingular and highly ill-conditioned, 
we shall prove that a stable numerical (affine subspace) solution may exist 
and contain an accurate approximation to the exact solution.
~For practical computation, efficient and robust algorithms already exist for 
general numerical solutions in affine Grassmannians.
~Regularization algorithms such as the Tikhonov method and
truncated SVD~\cite[\S 5.5.4]{golub-vanloan4}\cite{Han87} produce the 
accurate vector component and numerical rank-revealing algorithms%
~\cite{BarErbSla,utvtool}\cite[\S 5.4.6]{golub-vanloan4}%
\cite{lee-li-zeng,LeeLiZeng09,li-zeng-03,stew_utv}
provide the numerical kernel as the remaining component.

For the continuity of presentation, lemmas and long proofs are listed 
in the appendix.
~Additional computating results and software demonstration are given in the
supplementary material.

\section{Preliminaries}

Column vectors are denoted by boldface lower case letters such as ~$\bdb$,
~$\bdx$, ~$\bdy$ ~etc with ~$\bdo$ ~being a zero vector whose dimension can
be derived from the context.
~The vector space of ~$n$-dimensional complex column vectors is denoted by
~$\C^n$.
~The vector space of ~$m\times n$ ~matrices with complex entries is denoted
by ~$\C^{m\times n}$.
~Matrices are denoted by upper case letters such as ~$A$, ~$B$, ~$X$, ~etc
~with ~$O$ ~and ~$I$ ~denote a zero matrix and an identity matrix respectively.
~The range, kernel, rank and Hermitian transpose of a matrix ~$A$ ~are
denoted by ~$\cR(A)$, ~$\cK(A)$, ~$\rank{A}$ ~and ~$A^\h$ ~respectively.
~In this paper, we consider general ~$m\times n$ ~linear systems
in the form of ~$A\,\bdx\,=\,\bdb$ ~and we say the system is 
{\em singular} when ~$\rank{A}\,<\,n$ ~so that ~$\cK(A)\,\neq\,\{\bdo\}$,
~including non-square cases where ~$m\,<\,n$ ~or ~$m\,>\,n$.
~The system is {\em consistent} \,if ~$\bdb\,\in\,\cR(A)$.

For any matrix ~$A\,\in\,\C^{m\times n}$, 
~the ~$j$-th largest singular value of a matrix ~$A$ ~is denoted by
~$\sg_j(A)$.
~The {\em numerical rank}\, of a matrix ~$A$ ~{\em within an error tolerance}
~$\theta\,>\,0$ ~is defined as
\[ \ranka{\theta}{A} ~~:=~~ \min_{\|B-A\|_2<\theta}\rank{B}
~~\equiv~~ \max_{\sg_j(A)>\theta} j
\]
assuming ~$\theta$ ~does not equal to any singular value of ~$A$.
~Let ~$U\,\Sigma\,V^\h$ ~be the singular value decomposition of ~$A$
~where ~$U\,=\,[\bdu_1,\,\cdots,\,\bdu_m]$ ~and 
~$V\,=\,[\bdv_1,\,\cdots,\,\bdv_n]$.
~If ~$\ranka{\theta}{A} \,=\,r$ ~within ~$\theta$, ~then the 
{\em ~$\theta$-projection} ~$A_\theta$ ~of ~$A$ ~is defined as
\[ A_\theta ~~:=~~ \sg_1(A)\,\bdu_1\,\bdv_1^\h + \cdots +
\sg_r(A)\,\bdu_r\,\bdv_r^\h ~~=~~
\sum_{\sg_j(A)>\theta} \,\sg_j(A)\,\bdu_j\,\bdv_j^\h.
\]
In this case, the {\em numerical kernel}\, of 
~$A$ ~within ~$\theta$ ~is ~$\cK(A_\theta)\,=\,
\spn\{\bdv_{r+1},\,\ldots,\,\bdv_n\}$
~where ~$\spn\{\ldots\}$ ~denotes the vector space spanned by vectors in
the list.
~The entities ~$\ranka{\theta}{A}$, ~$A_\theta$ ~and ~$\cK(A_\theta)$ ~are 
undefined if ~$\theta$ ~is a singular value of ~$A$.
~The {\em Moore-Penrose inverse} ~of ~$A$, ~denoted by ~$A^\dagger$,
~is the unique matrix satisfying the Moore-Penrose conditions
~$A\,A^\dagger\,A\,=\, A$, ~$A^\dagger\,A\,A^\dagger\,=\, A^\dagger$, 
~$(A\,A^\dagger)^\h \,=\, A\,A^\dagger$ ~and
~$(A^\dagger\,A)^\h \,=\, A^\dagger\,A$.
~Using the singular value decomposition as above and assuming 
~$\rank{A}\,=\,r$, ~the identity~\cite[\S 5.5.2]{golub-vanloan4}
\[  A^\dagger ~~\equiv~~ 
\frac{1}{\sg_1(A)}\,\bdv_1\,\bdu_1^\h + \cdots +
\frac{1}{\sg_r(A)}\,\bdv_r\,\bdu_r^\h ~~=~~
\sum_{\sg_j(A)>0} \,\frac{1}{\sg_j(A)}\,\bdv_j\,\bdu_j^\h
\]
holds and ~$X\,=\,A^\dagger$ ~is the mininum Frobenius norm matrix ~such that 
~$A\,X$ ~and ~$X A$ ~are orthogonal projections from ~$\C^m$ ~onto ~$\cR(A)$ 
~and from ~$\C^n$ ~onto ~$\cR(A^\h)$ ~respectively.
~We shall frequently use ~$\big\|A^\dagger\big\|_2^{-1}$ ~as an alternative
notation for the smallest positive singular value ~$\sg_r(A)$ ~of ~$A$ 
~with rank ~$r$.

The set of \,$k$-dimensional subspaces of ~$\C^n$ ~is
called the {\em Grassmannian}~\cite{EdeAriSmi}\cite[page 52]{GraBook} 
of index ~$k$ ~of ~$\C^n$ ~denoted by ~$\cG_k(\C^n)$.
~For any ~${\cal P},\,{\cal Q}\,\in\,\cG_k(\C^n)$, ~let 
~$P,\,Q,\,\in\,\C^{n\times k}$ ~be matrices whose columns form 
orthonormal bases for ~${\cal P}$, ~${\cal Q}$ ~respectively while
~$\hat{P},\,\hat{Q}\, \in\,\C^{n\times(n-k)}$ ~such that 
~$[P,\hat{P}]^\h\,[P,\hat{P}]\,=\,[Q,\hat{Q}]^\h\,[Q,\hat{Q}]\,=\,I$.
~The Grassmannian ~$\cG_k(\C^n)$ ~is a metric space with the distance
\cite[\S 2.5.3]{golub-vanloan4}
\[  \dist{{\cal P},{\cal Q}} ~~:=~~ \big\|P\,P^\h-Q\,Q^\h\big\|_2
~~\equiv~~ \big\|P^\h\,\hat{Q}\|_2
~~\equiv~~ \big\|Q^\h\,\hat{P}\|_2.
\]
The set of ~$k$-dimensional affine 
subspaces of ~$\C^n$ ~is called the {\em affine Grassmannian}%
~\cite[\S 7.1]{AGraBook}\cite{naag,gras}
of index ~$k$ ~of ~$\C^n$ ~denoted by 
\[ \cA_k(\C^n) ~~:=~~ \big\{\bdu+\cV\,\subset\,\C^n ~\big|~ 
\bdu\,\in\,\C^n, ~\cV\,\in\,\cG_k(\C^n)\big\}.
\]
Here, for any vector  ~$\bdu\,\in\,\C^n$ ~and subspace 
~$\cV\,\in\,\cG_k(\C^n)$, ~the 
{\em affine subspace} 
\[ \bdu+\cV ~~:=~~ \{\bdu+\bdv\,\in\C^n ~|~ \bdv\,\in\,\cV\} 
\]
can be written as 
~$\hat\bdu+\cV$ ~with a unique ~$\hat\bdu\,\in\,\cV^\perp\cap
(\bdu+\cV)$ ~of the minimum norm where ~$(\cdot)^\perp$ ~denotes the 
unitary complement of any subspace ~$(\cdot)$.
~The metric
\begin{align}
\lefteqn{\dist{\bdu_1+\cV_1,\,\bdu_2+\cV_2}~~:=~~}  \nonumber \\
&
\label{afdis}
~~~~~~~\max_{\hat\bdu_j\in\cV_j^\perp\cap(\bdu_j+\cV_j),\,j=1,2}
\big\{\|\hat\bdu_1-\hat\bdu_2\|_2, ~\dist{\cV_1,\,\cV_2}\big\}
\end{align}
for every ~$\bdu_1+\cV_1,\,\bdu_2+\cV_2\,\in\,\cA_k(\C^n)$ 
~is a distance in ~$\cA_k(\C^n)$.
~For every ~$(A,\bdb)\,\in\,\C^{m\times n}\times\C^m$, ~denote
the set of vector solutions to the system ~$A\,\bdx\,=\,\bdb$ ~by
\[ \sol(A,\bdb) ~~:=~~ \{\bdu\,\in\,\C^n ~|~ A\,\bdu\,=\,\bdb\}.
\]
For ~$r\,=\,\rank{A}$, ~the set ~$\sol(A,\bdb)$ ~as {\em the} solution
of ~$A\,\bdx\,=\,\bdb$
~uniquely exists as either ~$\emptyset$ ~or
an element in the affine Grassmannian ~$\cA_{n-r}(\C^n)$. 
~The {\em dimension} of  ~$\sol(A,\bdb)$ ~is either ~$n-r$ ~if 
it is in ~$\cA_{n-r}(\C^n)$ ~or ~$-1$ ~if it is empty~\cite[page 6]{TopDim}.
~We define ~$\dist{\emptyset,\emptyset}\,=\,0$ ~so that the deviation of 
solutions can be measured if and only if they are of the same dimension.

The {\em condition number} of a square matrix ~$A$
~in the context of solving a linear system ~$A\,\bdx\,=\,\bdb$ ~is 
well-known to be ~$\kappa(A)\,=\,\|A\|_2\,\|A^{-1}\|_2$ ~with a convention 
~$\kappa(A)\,=\,\infty$ ~when ~$A$ ~is singular 
\cite[p. 87]{golub-vanloan4}.
~This condition number is based on the attainable error of the solution 
as an individual vector. 
~The infinity convention can be justified
by ~$\limsup_{G\rightarrow A}\|G\|_2\,\|G^\dagger\|_2 \,=\, \infty$
~when ~$A$ ~is singular
~and by the interpretation as the reciprocal of the distance to the 
singularity~\cite[Theorem 6.5]{Higham2}.
~For a rectangular ~matrix ~$A$, ~it is natural to generalize the 
condition number as ~$\kappa(A)\,=\,\|A\|_2\,\|A^\dagger\|_2$ 
~(see, e.g. \cite[p. 382]{Higham2}).
~We shall make arguments from multiple perspectives that the infinity 
convention may be unnecessary even if ~$A$ ~is square and singular.

It is easy to see that ~$\kappa(A)\,=\,\|A\|_2\,\|A^\dagger\|_2$
~is discontinuous at any rank deficient matrix ~$A$ ~and can not be 
approximated from empirical data ~$\tilde{A}$ ~since 
~$\kappa(\tilde{A})\,=\,\|\tilde{A}\|_2\,\|\tilde{A}^\dagger\|_2$
~can be arbitrarily large when ~$\|\Dl A\|_2\,=\,\|\tilde{A}-A\|_2$ ~is small.
~For any error tolerance ~$\theta$ ~with
~$0\,<\,\theta\,<\,\|A^\dagger\|_2^{-1}$, ~however, 
the asymptotic bound
\[  \kappa(A) - 2\,\|\Dl A\|_2 + O(\|\Dl A\|_2^2) ~~\le~~
\kappa(\tilde{A}_\theta) 
~~\le~~ \kappa(A) + 2\,\|\Dl A\|_2 + O(\|\Dl A\|_2^2)
\]
follows from \cite[Corollary 8.6.2]{golub-vanloan4}
when the data matrix ~$\tilde{A}$ ~is sufficiently accurate so that
~$\|\Dl A\|_2 \,<\,\|A^\dagger\|_2^{-1}-\theta$.
~Assuming an error bound ~$\bt\,>\,\|\Dl A\|_2$ ~is attainable
and is sufficiently tight so that ~$\bt\,<\,\|A^\dagger\|_2^{-1}-\|\Dl A\|_2$,
~the condition number ~$\kappa(\tilde{A}_\theta)$ ~of the ~$\theta$-projection
~$\tilde{A}_\theta$ ~of the data matrix ~$\tilde{A}$, 
~not ~$\kappa(\tilde{A})$,
~is an approximation to the underlying condition number 
~$\kappa(A)\,=\,\|A\|_2\,\|A^\dagger\|_2\,<\,\infty$.

\section{Models of singular linear systems}\label{s:mod}

We shall elaborate some case studies to show that
solving singular linear systems is not only unavoidable in scientific
computing, but also crucial in many applications.
~It may even be beneficial for the systems to be singular.
~Moreover, singular linear systems are often not known with exact matrices
and right-hand side vectors in practical computation, and need to be solved 
from empirical data.  

\begin{example}[Multiplicity of a singular solution to a nonlinear system]
\label{e:mul} \em
~For a system of nonlinear equations in the form of ~$\bdf(\bdx)\,=\,\bdo$
~where ~$\bdf \,=\,(f_1,\ldots,f_m)$ ~and
~$f_j~:~\C^n\,\rightarrow\,\C$ ~is an analytic function for
~$j\,=\,1,\ldots,m$, ~a zero ~$\bdx_*$ ~of ~$\bdf$ ~is multiple if the 
Jacobian of ~$\bdf$ ~at ~$\bdx_*$ ~is rank-deficient.
~At such a multiple ~$\bdx_*$ ~there is a vector space called the 
dual space ~$\cD_{\bdf,\bdx_*}$ ~that forms the multiplicity structure of the
zero and the dimension of ~$\cD_{\bdf,\bdx_*}$ ~is the multiplicity.
~The multiplicity structure can be determined by solving a
sequence of homogeneous linear systems 
\begin{equation}\label{sxc}  S_\al(\bdx_*)\,\bdc ~~=~~ \bdo
~~~~\mbox{for}~~~~ \al\,=\,1,2,\ldots
\end{equation}
where ~$S_\al(\bdx_*)$ ~is the Macaulay matrix whose entries are derivatives 
of ~$f_j$'s of orders up to ~$\al$ ~evaluated at ~$\bdx_*$.
~The solution ~$\sol\big(S_\al(\bdx_*),\bdo\big)$ ~of (\ref{sxc}) in a proper 
Grassmannian is isomorphic to the desired dual space ~$\cD_{\bdf,\bdx_*}$ 
~when ~$\al$ ~reaches the so-called depth.
~See, e.g., \cite{DLZ} for detailed elaborations and the supplemenary 
material for a computing demo.
~The exact Macaulay matrix ~$S_\al(\bdx_*)$ ~is almost never available since
~$\bdx_*$ ~is generally known approximately through a certain
~$\tilde\bdx\,\approx\,\bdx_*$ ~within an error bound.
~The model is to solve the singular system (\ref{sxc}) for 
the solution in a Grassmannian rather than individual vectors from empirical 
data matrix ~$S_\al(\tilde\bdx)\,\approx\,S_\al(\bdx_*)$.
\end{example}

\begin{example}[Sylvester equation] \label{e:syl} \em 
~This is an application arising in control and system theory\cite{AvrLas}
in the form of the Sylvester matrix equation
~$A(t)\,X+X\,B(t) \,=\, C(t)$
where ~$A(t)$, ~$B(t)$ ~and ~$C(t)$ ~are matrices depending on a 
parameter ~$t$.
~The system may inevitably become singular
when the parameter ~$t$ ~varies continuously and passes through a certain 
~$t_*$ ~whose value may only be obtained approximately.
~The following illustrative example is slightly modified from \cite{AvrLas}
(c.f. supplementary material). 
~Let
\begin{equation}\label{abct}  
A(t) ~=~ \mbox{\scriptsize $
\left[\begin{array}{rr} 1 & -1 \\ 1 & -1 \end{array}\right]$},
~~~
B(t) ~=~ \mbox{\scriptsize $
\left[\begin{array}{cc} -\frac{5}{3} + t & 1 \\ -1 & 
-\frac{1}{3} +2\,t \end{array}\right]$},
~~~\mbox{and}~~~
C(t) ~=~ \mbox{\scriptsize $
\left[\begin{array}{rr} 1 & 0 \\ 2 & -1 \end{array}\right]$}.
\end{equation}
When ~$t$ ~varies continuously, ~the system
becomes singular but still consistent when ~$t$ ~hits the value 
~$t_*\,=\,\frac{2}{3}$ ~with the general solution
\begin{equation}\label{751}  X_* ~~=~~ \mbox{$\frac{1}{4}$}\,
\mbox{\scriptsize $
\left[\begin{array}{rr} 1 & -1 \\ -3 & -1 \end{array}\right]$} + 
\al_1 \mbox{\scriptsize $
\left[\begin{array}{rr} 1 & 0 \\ 0 & 1 \end{array}\right]$} + 
\al_2 \mbox{\scriptsize $
\left[\begin{array}{rr} -1 & 1 \\ -1 & 1 \end{array}\right]$},
~~~~\al_1,\,\al_2\,\in\,\C.
\end{equation}
Suppose we know ~$\tilde{t}\,\approx\,0.6666$ ~with an
error bound ~$0.0001$.
~Can we find a numerical solution ~$\tilde{X}$ ~of the perturbed system at the 
parameter value ~$t\,=\,\tilde{t}$ ~approximating ~$X_*$ ~in (\ref{751}) 
of the underlying system at ~$t\,=\,t_*$ ~with an accuracy 
~$\|\tilde{X}-X_*\|_2$ ~roughly 0.0001?
\end{example}

\begin{example}[B\'ezout coefficients]\label{e:uf}   \em
~For polynomials ~$f_1,\,\ldots,\,f_n$, ~with a greatest common 
divisor ~$g$, ~there exist polynomials ~$u_1,\,\ldots,\,u_n$, ~known 
as the B\'ezout coefficients (see e.g. \cite[\S 1.3]{mora1}\cite{WilYu}), 
such that the B\'ezout identity
\begin{equation}\label{uf123}
 u_1\,f_1+\cdots+u_n\,f_n  ~~=~~ g
\end{equation}
holds.
~Solving the linear equation (\ref{uf123}) for the B\'ezout coefficients 
appears in many applications such as computing the Smith normal form
in linear control theory \cite{WilYu},
and the systems are often singular for ~$n\,\ge\,3$. 
~Denote ~$\mP_k$ ~as the vector space of polynomials with degrees up to ~$k$.
~For instance, ~let ~$f_1,\,f_2,\,f_3$ ~be polynomials of degrees, 
say ~$4,\,7,\,6$ ~with degree of ~$g$, ~say ~$2$, ~the equation
(\ref{uf123}) for ~$(u_1,\,u_2,\,u_3)\,\in\,\mP_3\times\mP_1\times\mP_2$ ~is
consistent and rank-deficient by 2. 
~The rank-deficiency is, in fact, a blessing in turning the general solution
\[ (u_1,\,u_2,\,u_3) ~=~ (u_{01},\,u_{02},\,u_{03})+
t_1\,(u_{11},\,u_{12},\,u_{13})+t_2\,(u_{21},\,u_{22},\,u_{23})
\]
into an invertible transformation 
\begin{equation}\label{ufg}
\mbox{\scriptsize $\left[\begin{array}{ccc}
u_{01} & u_{02} & u_{03} \\ 
u_{11} & u_{12} & u_{13} \\ 
u_{21} & u_{22} & u_{23} 
\end{array}\right]\,
\left[\begin{array}{c} f_{1} \\ f_{2} \\ f_{3} \end{array}\right]\,
$}
~~=~~
\mbox{\scriptsize $\left[\begin{array}{c} g \\ 0 \\ 0 \end{array}\right]$}.
\end{equation}
The exact coefficients of polynomial parameters 
~$f_1,\,\ldots,\,f_n$ ~and ~$g$ ~may be unknown beyond their empirical data, 
say
\begin{eqnarray*}
\tilde{f}_1 &~~=~~& 
\mbox{\scriptsize 2.5714} + \mbox{\scriptsize 3.8571}\,x - 
\mbox{\scriptsize 3}\,x^2 - \mbox{\scriptsize 6.4286}\,x^3 - 
\mbox{\scriptsize 2.1429}\,x^4 \\
\tilde{f}_2 &~~=~~& 
\mbox{\scriptsize -1.7143} - \mbox{\scriptsize 1.7143}\,x + 
\mbox{\scriptsize 0.4286}\,x^2 + \mbox{\scriptsize 0.4286}\,x^3 - 
\mbox{\scriptsize 3.4286}\,x^5 - \mbox{\scriptsize 5.1429}\,x^6 - 
\mbox{\scriptsize 1.7143}\,x^7 \\
\tilde{f}_3 &~~=~~& 
\mbox{\scriptsize 0.8571} + \mbox{\scriptsize 1.2857}\,x + 
\mbox{\scriptsize 2.1429}\,x^2 + \mbox{\scriptsize 2.5714}\,x^3 + 
\mbox{\scriptsize 3.4286}\,x^4 + \mbox{\scriptsize 3.8571}\,x^5 + 
\mbox{\scriptsize 1.2857}\,x^6 \\
\tilde{g} &~~=~~& \mbox{\scriptsize 4.6667} + \mbox{\scriptsize 7}\,x + 
\mbox{\scriptsize 2.3333}\,x^2
\end{eqnarray*}
with coefficientwise error bound ~$\eps\,=\,0.5\times 10^{-4}$.
~Can we accurately calculate the general solution for 
~$(u_1,\,u_2,\,u_3) \,\in\,\mP_3\times\mP_1\times\mP_2$ ~of the 
equation (\ref{uf123}) using the imperfect data 
~$\tilde{f}_1,\,\tilde{f}_2,\,\tilde{f}_3$ ~and ~$\tilde{g}$ ~within an error
in the same order of the data? 
~A computation/software demo for this example is given in the supplementary
material.
~The matrix-vector representation ~$\tilde{A}\,\bdx\,=\,\tilde\bdb$ ~of the 
equation (\ref{uf123}) in the given data is
\begin{equation} \label{54646}
\begin{array}{l}
\mbox{\tiny $\left[\!\begin{array}{rrrrrrrrr}
 2.5714 & 0 & 0 & 0 & -1.7143 & 0 &  0.8571 & 0 & 0 \\
 3.8571 &  2.5714 & 0 & 0 & -1.7143 & -1.7143 &  1.2857 &  0.8571 & 0 \\
-3.0000 &  3.8571 &  2.5714 & 0 &  0.4286 & -1.7143 &  2.1429 &  1.2857 &  0.8571 \\
-6.4286 & -3.0000 &  3.8571 &  2.5714 &  0.4286 &  0.4286 &  2.5714 &  2.1429 &  1.2857 \\
-2.1429 & -6.4286 & -3.0000 &  3.8571 & 0 &  0.4286 &  3.4286 &  2.5714 &  2.1429 \\
 0 & -2.1429 & -6.4286 & -3.0000 & -3.4286 & 0 &  3.8571 &  3.4286 &  2.5714 \\
 0 & 0 & -2.1429 & -6.4286 & -5.1429 & -3.4286 &  1.2857 &  3.8571 &  3.4286 \\
 0 & 0 & 0 & -2.1429 & -1.7143 & -5.1429 & 0 &  1.2857 &  3.8571 \\
 0 & 0 & 0 & 0 & 0 & -1.7143 & 0 & 0 &  1.2857
\end{array}\!\right]$}
\mbox{\tiny $\left[\!\begin{array}{c}
x_{1} \\ x_{2} \\ x_{3} \\ x_{4} \\ x_{5} \\ x_{6} \\ x_{7} \\ x_{8} \\ x_{9} 
\end{array}\!\right]$} \\
~=~ 
\mbox{\tiny $\left[\!\begin{array}{rrrrrrrrr}
 4.6667 & 7.0000 & 2.3333 & 0 & 0 & 0 & 0 & 0 & 0 
\end{array}\!\right]^\top$}
\end{array}
\end{equation}
with respect to monomial bases, and the condition number 
~$\kappa(\tilde{A}) \,\gtrapprox\, 2.29\times 10^6$.
~The system (\ref{54646}) in the conventional sense is highly 
ill-conditioned since ~$\eps\,\kappa(\tilde{A})\,>\,1$.
\end{example}

Applications are abundant involving singular linear systems. 
~The output regulation problem arises in the 
application of neural networks \cite{LiuHuang17}
for finding the matrix pair ~$(X,U)$ ~satisfying the so-called regulator 
equations whose solutions are not necessarily unique.
~An illustrative example is as follows (c.f. supplementary material):
\begin{equation} \label{xu}
\left\{\begin{array}{rcl}
X\,
\mbox{\scriptsize $\left[\begin{array}{cc}1 & 1 \\ 0 & 1 \end{array}\right]$}
& ~~=~~ & 
\mbox{\scriptsize $\left[\begin{array}{ccc} 0 & 1 & 0 \\ 0 & 0 & 1 \\ 
2 & -1 & 0 \end{array}\right]$}\,X + 
\mbox{\scriptsize $\left[\begin{array}{c} 0 \\ 0 \\ 1 \end{array}\right]$}\,U + 
\mbox{\scriptsize $\left[\begin{array}{rr}2 & 1 \\ -1 & 1 \\ 0 & 0 \end{array}
\right]$} \\
\mbox{\scriptsize $\left[\begin{array}{cc} 0 & 0 \end{array}\right]$}
 & = & 
\mbox{\scriptsize $\left[\begin{array}{ccc} 1 & 0 & -1 
\end{array}\right]$}\,X + 
\mbox{\scriptsize $\left[\begin{array}{cc} -1 & 0 \end{array}\right]$}
\end{array}\right.
\end{equation}
where the unknowns ~$X$ ~and ~$U$ ~are matrices.
~The system is rank deficient by one.
~Furthermore, the matrix parameters 
are not known exactly but given by estimation.
~For a matrix ~$A$ ~with a defective
eigenvalue ~$\la_*$ ~and an associated eigenvector ~$\bdz_*$, 
~a generalized eigenvector satisfies the singular system
~$(A-\la_*\,I)\,\bdx \,=\, \bdz_*$ ~for ~$\bdx\,\in\,\C^n$.
~The value of ~$\la_*$ ~and ~$\bdz_*$ ~generally can
only be known approximately.
~The problem is to solve the underlying system by solving 
~$(\tilde{A}-\tilde\la\,I)\,\bdx\,=\,\tilde\bdz$ ~from
the data ~$\tilde{A}\,\approx A$, ~$\tilde\la\,\approx\,\la_*$
~and ~$\tilde\bdz\,\approx\,\bdz_*$.
~More applications include solving the singular homogeneous linear systems 
of Ruppert matrices in numerical factorization of polynomials 
\cite{gao_fac,WuZeng}, numerical elimination of polynomial variables
\cite{ZengElim}, etc.
~The generalized Lyapunov equation
~$E^\h\,A\,X + A^\h\,X\,E \,=\, -G$
~with given matrices ~$A$, ~$E$ ~and ~$G$ ~is singular when ~$E$ ~is
rank-deficient \cite{sty02}.
~A singular linear system that models the atmospheric path delay and the water
vapor constant estimation is given in \cite{SaqKar}.
~Linear systems derived from discretizing the Fredholm and Volterra
integral equations can be considered empirical data of singular systems in the 
presense of annihilators \cite[\S 2.4 and page 83]{HansInv}
(c.f. an example in supplementary material).

\section{Homogeneous systems with empirical data}\label{s:hs}

A problem is {\em well-posed}\, if its solution satisfies existence, 
uniqueness and Lipschitz continuity with respect to the data 
or, otherwise, it is an {\em ill-posed problem}.
~For an ~$m\times n$ ~singular homogeneous linear system 
~$A\,\bdx\,=\,\bdo$, ~the problem
\begin{equation}\label{pax0}
\mbox{Solve ~$A\,\bdx\,=\,\bdo$ ~for a single-vector solution 
~$\bdx$ ~in ~$\C^n$}
\end{equation}
is obviously ill-posed as its solutions are not unique.
~However, the problem (\ref{pax0}) is not precisely the problem to be solved 
in standard linear algebra where all the solutions are in question.
~There is a unique solution to the problem
\begin{equation}\label{psa0}
\mbox{Solve ~$A\,\bdx\,=\,\bdo$ ~for the solution ~$\sol(A,\bdo)$ ~in
the Grassmannian ~$\cG_{n-r}(\,\C^n)$} 
\end{equation}
where ~$r\,=\,\rank{A}$.
%
~The problem may become somewhat confounding when the exact ~$A$ ~is 
unknown but given through empirical data in 
~$\tilde{A}$ ~as illustrated in Example~\ref{e:mul}.
~What really is at stake is a nontrivial solution ~$\sol(A,\bdo)\,\equiv\,
\cK(A)$ ~in the
Grassmannian ~$\cG_{n-r}(\C^n)$ ~but the data system 
~$\tilde{A}\,\bdx\,=\,\bdo$ 
~is almost always nonsingular with ~$\sol(\tilde{A},\bdo)\,=\,\{\bdo\}\,\in\,
\cG_0(\C^n)$ ~when ~$m\,\ge\,n$.
%
%
~The condition number
~$\kappa(\tilde{A})\,=\,O(\|A-\tilde{A}\|_2^{-1})$ ~can be huge as well 
if ~$r\,<\,\min\{m,\,n\}$.
~The very problem of solving a homogeneous linear system from empirical data
needs clarification.
%

\begin{problem}[Numerical Solution of a Homogeneous Linear System]\label{p:rr}
~Let ~$\tilde{A}$ ~be an ~$m\times n$ ~matrix serving as empirical 
data for an underlying homogeneous system ~$A\,\bdx\,=\,\bdo$ ~where 
entries of ~$A$ ~may or may not be known exactly. 
~Identify the rank ~$r$ ~of ~$A$ ~using ~$\tilde{A}$ ~and find a numerical 
solution of ~$\tilde{A}\,\bdx\,=\,\bdo$ ~ in the Grassmannian 
~$\cG_{n-r}(\C^n)$ ~in the form of an orthonormal basis 
~$\{\bdz_1,\ldots,\bdz_{n-r}\}$ ~so that 
\begin{equation}\label{nk}
\dist{\spn\{\bdz_1,\ldots,\bdz_{n-r}\},~\sol(A,\bdo)} ~~=~~
O\left(\mbox{$\frac{\|A-\tilde{A}\|_2}{\|A\|_2}$}\right).
\end{equation}
\end{problem}

From Wedin's perturbation analysis \cite{wedin}, the numerical kernel 
~$\cK(\tilde{A}_\theta)$ ~within a proper error tolerance ~$\theta\,>\,0$ 
~is an approximation to ~$\cK(A)\,=\,\sol(A,\bdo)$ ~in ~$\cG_{n-r}(\C^n)$
(c.f. Lemma~\ref{l:d} in \S\ref{s:lem} in appendix).
~For every ~$G\,\in\,\C^{m\times n}$, ~we define
\[ \sol_\theta(G,\bdo) ~~:=~~ \cK(G_\theta) ~~\equiv~~ 
\sol(G_\theta,\bdo)
\]
as the {\em numerical solution} of the homogeneous system ~$G\,\bdx\,=\,\bdo$
~in the Grassmannian ~$\cG_{n-r}(\C^n)$ ~{\em within an error tolerance} 
~$\theta$ ~where ~$r\,=\,\ranka{\theta}{G}$.
~Numerical methods for computing ~$\cK(G_\theta)$ ~as ~$\sol_\theta(G,\bdo)$ 
~are well-established, including the singular value decomposition and
other numerical rank-revealing methods 
(see, e.g. \cite{HansenBook,li-zeng-03}).
~The following theorem summarizes the properties of the 
numerical solution as a generalization of the exact solution to the 
homogeneous system and as a well-posed computing problem that solves the 
underlying system in Problem~\ref{p:rr}.
~The essence and underlying substance of Theorem~\ref{t:rr} are based on 
Wedin \cite{wedin}. 

\begin{theorem}\label{t:rr}
~Let ~$A\,\in\,\C^{m\times n}$. 
~The following properties hold for the numerical solution 
of a homogeneous system.
\begin{itemize}\parskip0mm
\item[\em (i)] ~{\em The exact solution is a special case of the
numerical solution:} 
\[  \sol(A,\bdo) ~~\equiv~~ \sol_\theta(A,\bdo) 
~~~~\mbox{for all}~~~~\theta\,\in\,\big(0,\|A^\dagger\|_2^{-1}\big). 
\]
\item[\em (ii)] ~{\em Computing the numerical solution is a well-posed 
problem:}
~If ~$\sol_\theta(A,\bdo)$ ~is well-defined within ~$\theta\,>\,0$, 
~then  ~$\sol_\theta(A+\Dl A,\bdo)$ ~uniquely exists in the same Grassmannian
as ~$\sol_\theta(A,\bdo)$ ~and enjoys Lipschitz continuity with
\begin{align} \lefteqn{
\dist{\sol_\theta(A+\Dl A,\,\bdo),\,\sol_\theta(A,\bdo)}} \nonumber \\
& ~~\le~~
\frac{\|A_\theta\|_2\,\big\|A_\theta^\dagger\|_2}{
1-\big\|A_\theta^\dagger\|_2\,(\|A-A_\theta\|_2+\|\Dl A\|_2)}\, 
\frac{\|\Dl A\|_2}{\|A\|_2}  \label{dsta0}
\end{align}
for all ~$\Dl A$ ~with sufficiently small ~$\|\Dl A\|_2$ ~satisfying
\[ \|\Dl A\|_2 ~~\le~~
\min\left\{\mbox{$\frac{1}{2}$}\,\left(\|A_\theta^\dagger\|_2^{-1}-
\|A-A_\theta\|_2\right),\,\theta-\|A-A_\theta\|_2,\,
\|A_\theta^\dagger\|_2^{-1}-\theta\right\}.
\]
\item[\em (iii)] 
~{\em A homogeneous system can be solved from empirical data with an 
accuracy in the same order as the data:} 
~For any ~$A+\Dl A$ ~serving as empirical data of ~$A$ ~with 
~$\|\Dl A\|_2\,<\, \frac{1}{2}\,\|A^\dagger\|_2^{-1}$, 
~there exist ~$\mu,\,\eta\,>\,0$ ~with
\begin{equation}\label{tbd0}
\mu ~~\le~~ \|\Dl A\|_2~~<~~ \|A^\dagger\|_2^{-1}-\|\Dl A\|_2
~~\le~~ \eta
\end{equation}
such that the numerical solution ~$\sol_\theta(A+\Dl A,\,\bdo)$ 
~within any error tolerance ~$\theta\,\in\,(\mu,\,\eta)$
~is in the same Grassmannian as the exact solution ~$\sol(A,\bdo)$ 
~and
\begin{equation}\label{dsa0}
\dist{\sol_\theta(A+\Dl A,\,\bdo),\,\sol(A,\bdo)}  ~~\le~~
\frac{\|A\|_2\,\big\|A^\dagger\|_2}{1-\big\|A^\dagger\|_2\,
\|\Dl A\|_2}\, \frac{\|\Dl A\|_2}{\|A\|_2}.
\end{equation}
\end{itemize}
\end{theorem}

{\em Proof.}
~A straightforward verification from Wedin's error bound \cite{wedin} on 
singular subspaces (see Lemma~\ref{l:d} in Appendix~\ref{s:lem}) along with 
the identity
~$A_\theta\,\equiv\,A$ ~for ~$0\,<\theta\,<\,\|A^\dagger\|_2^{-1}
\,=\,\sg_r(A)$ ~where ~$r\,=\,\rank{A}$,
~$\mu\,=\,\sg_{r+1}(\tilde{A})\,\le\,\|\Dl A\|_2$ ~and
~$\eta\,=\,\sg_r(\tilde{A})\,\ge\,\sg_r(A)-\|\Dl A\|_2$. ~\Qed

By Theorem~\ref{t:rr}, Problem~\ref{p:rr} is solvable if
the data are sufficiently accurate and a tight error bound on data
is attainable, as asserted in the following corollary.

\begin{corollary}\label{c:rr}
~Let the matrices ~$A$ ~and ~$\tilde{A}$ ~be as in Problem~\ref{p:rr}. 
~Assume the data in ~$\tilde{A}$ ~are sufficiently accurate
such that ~$\|A-\tilde{A}\|_2\,<\,\frac{1}{2}\,\|A^\dagger\|_2^{-1}$.
~Further assume a data error bound ~$\beta\,>\,\|A-\tilde{A}\|_2$ ~is known
and is sufficiently tight so that 
~$\beta\,<\,\|A^\dagger\|_2^{-1}-\|A-\tilde{A}\|_2$.
~Then Problem~\ref{p:rr} is solvable by setting the error tolerance 
~$\theta\,=\,\beta$ ~and finding an orthonormal basis for the numerical
solution ~$\sol_\theta(\tilde{A},\bdo)\,=\,\cK(\tilde{A}_\theta)$ ~within 
~$\theta$.
\end{corollary}

{\em Proof.}
~A straightforward verification using Theorem~\ref{t:rr}.
~\Qed

The error tolerance ~$\theta$ ~in Theorem~\ref{t:rr} is
an operational parameter that needs to be set up for solving 
Problem~\ref{p:rr}.
~If we assume the underlying application allows the data error to a certain
extent, say ~$\|A-\tilde{A}\|_2\,<\,\hat\theta$, ~the data error bound ~$\bt$
~in Corollary~\ref{c:rr} ~is expected to be below ~$\hat\theta$.
~The inequality (\ref{tbd0}) ensures there is a window ~$(\mu,\,\eta)$ ~for 
setting the operational error tolerance ~$\theta$ ~at ~$\beta$ ~or slightly 
larger.
~Using the notation of Problem~\ref{p:rr}, it is reasonable to assume 
the data error bound ~$\beta$ ~on ~$\|A-\tilde{A}\|_2$ ~is known or can be 
estimated.
~The crucial criterion for operational purpose is 
to set ~$\theta$ ~at or slightly above ~$\|A-\tilde{A}\|_2$ ~according
to Theorem~\ref{t:rr}, part (iii). 
~The error tolerance ~$\theta$ ~should not exceed 
~$\|A^\dagger\|_2^{-1}-\|A-\tilde{A}\|_2$ ~whose exact value or estimation
is not needed if the data error bound ~$\beta$ ~is sufficiently tight.
~See the supplementary material for examples of setting error tolerances.

For a rank-$r$ ~matrix ~$A$, ~the sensitivity of solving 
~$A\,\bdx\,=\, \bdo$ ~for 
~$\sol_\theta(\tilde{A},\bdo)$ ~in the Grassmannian ~$\cG_{n-r}(\C^n)$ 
~from a perturbed data matrix ~$\tilde{A}$ ~is
\[  \big\|A\big\|_2\,\big\|A^\dagger\big\|_2
~~=~~ \frac{\sg_1(A)}{\sg_r(A)} ~~\approx~~
  \big\|\tilde{A}_\theta\big\|_2\,\big\|\tilde{A}_\theta^\dagger\big\|_2
\]
from (\ref{dsta0}) and (\ref{dsa0}), not infinity or ~$\kappa(\tilde{A})$.
~The convention ~$\kappa(A)\,=\,\infty$ ~for the square singular case
and ~$\kappa(\tilde{A})\,=\,\|\tilde{A}\|_2\,\|\tilde{A}^\dagger\|_2$ ~may 
overestimate the sensitivity substantially.
~Problem~\ref{p:rr} may not be solvable if the data error 
is large beyond, say ~$\frac{1}{2}\,\|A^\dagger\|_2^{-1}$,
~or may not be solved accurately if the data error bound is unknown or
the inherent sensitivity ~$\|A\|_2\,\|A^\dagger\|_2$ ~is high. 
%

For solving ~$A\,\bdx\,=\,\bdo$ ~with ~$A\,\in\,\C^{m\times n}$, 
~there are differences between cases of ~$m\,<\,n$ ~and ~$m\,\ge\,n$.
~The solution is of a positive dimension when ~$m\,<\,n$ ~regardless of 
perturbations and, if ~$\rank{A}\,=\,m$, ~the condition
~$\|A\|_2\,\|A^\dagger\|_2$ ~is continuous with respect to small perturbations.
~When ~$m\,\ge\,n$ ~and ~$\sol(A,\,\bdo)$ ~is nontrivial, however, the 
dimension of ~$\sol(A+\Dl A,\,\bdo)$ ~degrades to zero for almost all 
perturbations ~$\Dl A$ ~and the condition ~$\|A\|_2\,\|A^\dagger\|_2$ ~is
discontinuous.
~The assertions of Theorem~\ref{t:rr} ~remain the same either ~$m\,<n$
~or ~$m\,\ge\,n$.

\section{Sensitivity of a consistent singular system} \label{s:sls}

Solving a singular system for an individual vector solution is known to have an
unbounded sensitivity under arbitrary perturbations. 
~From a different perspective, the infinity condition number is not the 
sensitivity of the {\em singular system} if the singularity is not maintained. 
~There is an intrinsic stability in solving
~$A\,\bdx\,=\,\bdb$ ~when the rank and consistency are preserved. 
~This point of view is originated in \cite{kahan72} by Kahan who suggests the 
perceived hypersensitivity of multiple roots may be a ``misconception'' 
without maintaining the multiplicity.

A consistent ~$m\times n$ ~linear system ~$A\,\bdx\,=\,\bdb$ ~with 
~$\rank{A}\,=\,r$
~has a unique solution ~$\sol(A,\bdb) \,=\, \bdx_0+\cK(A)$
~in the affine Grassmannian ~$\cA_{n-r}(\C^n)$
~where ~$\bdx_0$ ~is any particular solution.
~The sensitivity of the linear system ~$A\,\bdx\,=\,\bdb$ ~can be based on 
the deviation of the solution ~$\sol(A,\bdb)$ ~in ~$\cA_{n-r}(\C^n)$ ~with 
respect to perturbations of ~$(A,\bdb)\,\in\,\C^{m\times n}\times\C^n$.
~From (\ref{afdis}), the difference between solutions of two consistent 
systems of the same rank can be measured by the metric (\ref{afdis}), namely
\begin{align}
\lefteqn{\dist{\sol(A,\bdb),\,\sol(B,\bdd)} ~=~}\nonumber \\
& ~~~~~~~~\max\big\{\|A^\dagger\,\bdb-B^\dagger\,\bdd\|_2, 
~\dist{\cK(A),\,\cK(B)}\big\}.    \label{lsdist}
\end{align}
Notice that the component 
~$\dist{\cK(A),\,\cK(B)}\,\le\,1$ ~in (\ref{lsdist}) but the other 
component ~$\|A^\dagger\,\bdb-B^\dagger\,\bdd\|_2$ 
~can be large or small. 
~One way to avoid an imbalance is to put a weight factor ~$\omega$ ~on the 
component ~$\|A^\dagger\,\bdb-B^\dagger\,\bdd\|_2$.
~We choose not to use weights for the sake of simplicity of
elaborations and for the reason that the weight ~$\omega$ ~can be used
to scale the linear system instead so that we can solve 
~$A\,(\omega\,\bdx)\,=\,\omega\,\bdb$ ~equivalently.
~For convenience, we adopt a specific norm
\begin{equation}\label{Abnorm}
 \|(A,\bdb)\| ~~:=~~ \sqrt{\|A\|_2^2+\|\bdb\|_2^2} 
\end{equation}
in the product space ~$\C^{m\times n}\times\C^n$.
~The theories in this paper can be adapted to other norms.

With these notations and metrics, the solution ~$\sol(A,\bdb)$ ~of a 
singular consistent ~$m\times n$ ~linear system ~$A\,\bdx\,=\,\bdb$ ~uniquely 
exists in the affine Grassmannian ~$\cA_{n-r}(\C^n)$ ~and the sensitivity is 
proportional to ~$\|A\|_2\,\|A^\dagger\|_2$ ~rather than infinity when the 
rank and consistency are preserved, as established in the following 
theorem.

\begin{theorem}\label{t:lsps}
~{\em The solution of a consistent linear system is Lipschitz 
continuous when the rank and consistency are preserved.}
~Let ~$A\,\in\,\C^{m\times n}$ ~and ~$\bdb\,\in\,\cR(A)$. 
~Assume the perturbation ~$(\Delta A,\,\Delta \bdb)$ ~is constrained
such that:
~$\rank{\tilde{A}}\,=\,\rank{A}$ ~and 
~$\tilde\bdb\,\in\,\cR(\tilde{A})$ ~where ~$\tilde{A}\,=\,A+\Dl A$
~and ~$\tilde\bdb\,=\,\bdb+\Dl\bdb$.
~Then
\begin{align} 
\lefteqn{\dist{\sol(\tilde{A},\,\tilde\bdb),\,\sol(A,\,\bdb)}} \nonumber \\
 & ~~\le~~  \|A\|_2\,\|A^\dagger\|_2\cdot
\frac{\sqrt{2\,\|\bdx_*\|_2^2+1}}{
\|A\|_2-\sqrt{2}\,\|A\|_2\,\|A^\dagger\|_2\, \|\Dl A\|_2
}\, \|(\Dl A,\,\Dl \bdb)\|  \label{afsest}
\end{align}
where ~$\bdx_*\,=\,A^\dagger\,\bdb$ ~whenever 
~$\sqrt{2}\,\|A^\dagger\|_2\,\|\Dl A\|_2\,<\,1$.
\end{theorem}

{\em Proof sketch.} 
~The kernel component of the distance in (\ref{afsest}) is bounded by
Wedin's error estimate \cite{wedin} (see Lemma~\ref{l:d} in
Appendix~\ref{s:lem}).
~Let ~$N$ ~be a matrix whose coluns form an orthonormal basis for 
~$\cK(A)$. 
~Then the mininum norm solution ~$\bdx_*$ ~is the unique least squares
solution of the system 
\[\left[\begin{array}{c} \mu\,N^\h \\ A \end{array}\right]\,\bdx ~~=~~
\left[\begin{array}{c} \bdo \\ \bdb \end{array}\right]
~~~~\mbox{for any}~~~~ \mu\,>\,0
\]
and the standard error bound \cite[Theorem 1.4.6]{bjorck96} applies.
~Detailed proof is in Appendix~\ref{s:prf}.
~\Qed

As a result of (\ref{afsest}), the intrinsic sensitivity of solving 
a singular system ~$A\,\bdx\,=\,\bdb$ ~for 
the general solution ~$\sol(A,\bdb)$ ~is a constant multiple of
\[  \|A\|_2\,\big\|A^\dagger\big\|_2 ~~=~~ \frac{\sg_1(A)}{\sg_r(A)}
~~<~~ \infty
\]
when the rank and consistency are preserved.

As a by-product of establishing Theorem~\ref{t:lsps}, the following corollary 
improves the standard normwise
error bound \cite[Theorem 5.6.1]{golub-vanloan4}
on the minimum norm solution of a full rank underdetermined linear system
by reducing a factor from 2 to ~$\sqrt{2}$.

\begin{corollary}\label{c:uds}
~Let ~$A\,\in\,\C^{m\times n}$ ~with ~$\rank{A}\,=\,m\,<\,n$ ~and
~$\bdb\,\in\,\C^m$. 
~If ~$\bdx_*$ ~and ~$\tilde\bdx$ ~are minimum norm solutions of the
underdetermined linear systems ~$A\,\bdx\,=\,\bdb$ ~and
~$(A+\Dl A)\,\bdx\,=\,\bdb+\Dl\bdb$ ~respectively with
~$\sqrt{2}\,\|A^\dagger\|_2\,\big\|\Dl A\big\|_2\,<\,1$, 
~then
\begin{align} \frac{\|\tilde\bdx-\bdx_*\|_2}{\|\bdx_*\|_2} &~~\le~~ 
\|A\|_2\,\|A^\dagger\|_2\,\left(\sqrt{2}\,
\frac{\|\Dl A\|_2}{\|A\|_2}+\frac{\|\Dl \bdb\|_2}{\|\bdb\|_2}
\right) + O(\|(\Dl A,\Dl \bdb)\|^2). \label{udms}
\end{align}
\end{corollary}

{\em Proof.} 
~The inequality (\ref{udms}) ~follows from (\ref{26428}) in 
Appendix~\ref{s:prf}.  ~\Qed

\begin{remark}\em
The subset of all rank-$r$ matrices is a complex analytic 
manifold in the topological space ~$\C^{m\times n}$ \cite{dem-edel}
with the topology derived from the Frobenius norm. 
~Similarly the subset
~$\cM^{m\times n}_r \,:=\, \big\{(A,\bdb)\in\C^{m\times n}\times\C^m
~\big|~ \rank{A}=r,~\bdb\in\cR(A)\big\}$
~is a complex analytic manifold in ~$\C^{m\times n}\times\C^m$.
~Although the problem of solving a singular linear system in general is 
ill-posed, Theorem~\ref{t:lsps} implies the problem of solving 
~$A\,\bdx\,=\,\bdb$ ~for ~$\sol(A,\bdb)$ ~in ~$\cA_{n-r}(\C^n)$ 
~is well-posed on the manifold ~$\cM^{m\times n}_r$.
\end{remark}

\section{The general numerical solution}\label{s:gns}

When a rank-deficient ~$m\times n$ ~linear system ~$A\,\bdx\,=\,\bdb$ ~is 
given through empirical data ~$(\tilde{A},\tilde\bdb)$, ~the perturbed matrix
~$\tilde{A}$ ~is almost always of full rank and highly ill-conditioned.
~Furthermore, the conventional single-vector solution of the data system 
~$\tilde{A}\,\bdx\,=\,\tilde\bdb$ ~is in ~$\C^n$ ~while the general solution
of the underlying system is in completely different space ~$\cA_{n-r}(\C^n)$. 
~What the problem precisely is and what the numerical solution really means
need to be clarified.

\begin{problem}[Numerical Solution of a Linear System]\label{p:lsgs}
~For given ~$\tilde{A}$ ~and ~$\tilde\bdb$ 
~serving as empirical data for an underlying linear system 
~$A\,\bdx\,=\,\bdb$ ~to be solved, ~find a numerical solution of 
~$\tilde{A}\,\bdx\,=\,\tilde\bdb$ 
~that can be identified as the exact solution ~$\sol(\hat{A},\hat\bdb)$ 
~of ~$\hat{A}\,\bdx\,=\,\hat\bdb$
~such that both the backward error and the forward error 
\begin{align}\label{lsber}
  \big\|(\tilde{A},\tilde\bdb)-(\hat{A},\,\hat\bdb)\big\| &~~=~~
O\big(\big\|(\tilde{A},\,\tilde\bdb)-(A,\,\bdb)\big\|\big) \\
\label{lsfer}  
\dist{\sol(\hat{A},\hat\bdb),\,\sol(A,\,\bdb)} &~~=~~
O\big(\big\| (\tilde{A},\,\tilde\bdb) - (A,\,\bdb) \big\|\big)
\end{align}
are in the same order of the data accuracy.
\end{problem}

The accuracy requirement (\ref{lsfer}) stipulates that both
~$\sol(A,\bdb)$ ~and ~$\sol(\hat{A},\,\hat\bdb)$ 
~are in the same affine Grassmannian or both empty.
~It is natural to choose ~$\hat{A}\,=\,\tilde{A}_\theta$ ~within a proper
~$\theta$ ~and ~$\hat\bdb\,=\,\tilde\bdb_\theta\,:=\,
\tilde{A}_\theta\,\tilde{A}_\theta^\dagger\,\tilde\bdb$ ~as the orthogonal 
projection of ~$\tilde\bdb$ ~onto ~$\cR(\tilde{A}_\theta)$.
~The solution ~$\sol(\tilde{A}_\theta,\tilde\bdb_\theta)$ ~is acceptable as 
the numerical 
solution of ~$\tilde{A}\,\bdx\,=\,\tilde\bdb$ ~if its backward error
is below the error tolerance or the empty set.

\begin{definition}[General Numerical Solution]\label{d:ngs}
~Let ~$G\,\in\,\C^{m\times n}$, ~$\bdd\,\in\,\C^m$ ~and ~$\theta\,>\,0$
be an error tolerance within which ~$\ranka{\theta}{G}$ ~is well defined.
~With respect to a norm ~$\vvv\cdot\vvv$ ~on
~$\C^{m\times n}\times\C^m$, ~the {\em general numerical solution}
\,of the linear system ~$G\,\bdx\,=\,\bdd$ ~{\em within} ~$\theta$ ~is 
defined as
\[ \sol_\theta(G,\bdd) ~~:=~~ \left\{\begin{array}{ccl}
\sol(G_\theta,\bdd_\theta)
& &
\mbox{if ~$\vvv  (G,\bdd)-(G_\theta,\bdd_\theta)\vvv\,<\,\theta$} \\
\emptyset && 
\mbox{if ~$\vvv(G,\bdd)-(G_\theta,\bdd_\theta)\vvv\,>\,\theta$} \\
\end{array}\right.
\]
where ~$G_\theta$ ~is the ~$\theta$-projection of ~$G$
and ~$\bdd_\theta\,=\,G_\theta\,G_\theta^\dagger\,\bdd$ ~is 
the orthogonal projection of ~$\bdd$ ~onto the range
~$\cR(G_\theta)$ ~of ~$G_\theta$. 
\end{definition}

The solution ~$\sol_\theta(G,\bdd)$ ~is undefined 
if ~$\theta$ ~equals to a singular value of ~$G$ ~or 
~$\theta\,=\,\vvv(G,\bdd)-(G_\theta,\bdd_\theta)\vvv$.
%
%
%
~We can now establish the following theorem on 
the general numerical solution.

\begin{theorem}\label{t:ftngs}
~At any ~$(A,\,\bdb)\,\in\,\C^{m\times n}\times\C^m$,
~the following properties of the general numerical solution hold 
with respect to the norm {\em (\ref{Abnorm})}.
\begin{itemize}\parskip0mm
\item[\em (i)] ~{\em An exact general solution is a special case of 
general numerical solution}: ~The identity 
~$\sol(A,\,\bdb) \,\equiv\, \sol_\theta(A,\,\bdb)$ ~holds 
for all ~$\theta \,<\,\|A^\dagger\|_2^{-1}$ ~if ~$\bdb\,\in\, \cR(A)$, 
~or ~$\theta\,<\,\min\left\{\|A^\dagger\|_2^{-1}, 
\,\|\bdb-A\,A^\dagger\,\bdb\|_2\right\}$ ~otherwise.
\item[\em (ii)] ~{\em Computing the general numerical solution is a well-posed
problem:} ~Assume ~$\sol_\theta(A,\bdb)$ ~is well-defined within a certain 
~$\theta\,>\,0$.
~There is a ~$\xi\,>\,0$ ~depending on ~$A,\,\bdb$ ~and ~$\theta$
~such that, for every ~$(\Dl A,\,\Dl\bdb)$ ~with a sufficiently small norm,
~there exists a unique ~$\sol_\theta(A+\Dl A,\,\bdb+\Dl\bdb)$ 
~satisfying the Lipschitz continuity
\begin{equation}
\dist{\sol_\theta(A+\Dl A,\,\bdb+\Dl\bdb),\,\sol_\theta(A,\bdb)} ~~\le~~
\xi\,\big\|(\Dl A,\,\Dl\bdb)\big\|.
\label{lip}
\end{equation}
\item[\em (iii)] ~~{\em A singular linear system can be solved from 
empirical data with an accuracy in the same order as the data:} 
\newline {\em (a)}
~Assume ~$\bdb\,\in\,\cR(A)$ ~and let ~$\bdx_*\,=\,A^\dagger\,\bdb$. 
~For any empirical data pair ~$(\tilde{A},\,\tilde\bdb)\,=\,
(A+\Dl A,\,\bdb+\Dl\bdb)$ ~satisfying
~$\big\|(\Dl A,\,\Dl\bdb)\big\| \,<\,
\big((\omega+1)\,\|A^\dagger\|_2\big)^{-1}$
~where 
~$\omega\,=\, \sqrt{4\,\|A^\dagger\|_2^2\,\|\bdb\|_2^2+2}$
~and for any error tolerance ~$\theta$ ~satisfying
\begin{equation}\label{tintv}
  \omega\,\big\|(\Dl A,\,\Dl\bdb)\big\| ~~<~~ \theta
~~<~~ \|A^\dagger\|_2^{-1}- \big\|(\Dl A,\,\Dl\bdb)\big\|,
\end{equation}
there exists a unique general numerical solution 
~$\sol_\theta(\tilde{A},\tilde\bdb)$ 
~with a backward error bound
~$\omega\,\big\|(\Dl A,\,\Dl\bdb)\big\|$ ~and a forward error bound
\begin{align}
& \dist{\sol_\theta(\tilde{A},\tilde\bdb),\, \sol(A,\bdb)} \nonumber \\
& ~~~~~\le~  
\|A\|_2\,\|A^\dagger\|_2\,\frac{
\sqrt{4\,\|\bdx_*\|_2^2+1}
}{
\|A\|_2 -\|A\|_2\,\|A^\dagger\|_2\,\|\Dl A\|_2}\,\big\|(\Dl A,\,\Dl\bdb)\big\|.
\label{ngserr}
\end{align}
{\em (b)} ~Assume ~$\sol(A,\bdb)\,=\,\emptyset$. 
~For any ~$\theta\,<\,
\min\left\{\frac{1}{2}\,\|A^\dagger\|_2^{-1},\, \|\bdb-A\,A^\dagger\,\bdb\|_2
\right\}$,
~there is a constant ~$\rho\,\in\,(0,\theta)$ ~such that
~$\sol_\theta(\tilde{A},\,\tilde\bdb)\,=\,\sol(A,\bdb)\,=\,\emptyset$ 
~at any empirical data 
pair ~$(\tilde{A},\,\tilde\bdb)$ ~satisfying
~$\big\|(\tilde{A},\,\tilde\bdb)-(A,\,\bdb)\big\| \,<\, \rho$.
\end{itemize}
\end{theorem}

{\em Proof sketch.} 
~The assertion (i) and the unique existence in the assertion (ii) directly 
follow from Definition~\ref{d:ngs}.
~The Lipschitz continuity (\ref{lip}) is a variation of 
the error estimate for the truncated SVD solution by 
Hansen~\cite[inequality (26a)]{Han87} as an extension of Wedin error 
analysis \cite{wedin73}.
~The bound on the minimum norm solution component of the distance in 
the inequality (\ref{ngserr}) follows from Hansen 
~\cite[inequality (27a)]{Han87} and the 
the bound on the numerical kernel is established by Wedin \cite{wedin}.
~Detailed proof is given in Appendix~\ref{s:prf}. \Qed

%
For Problem~\ref{p:lsgs}, assume the underlying linear system 
~$A\,\bdx\,=\,\bdb$ ~in Problem~\ref{p:lsgs} is known to be consistent 
in applications such as Example~\ref{e:uf}, the
solvability the system from empirical data ~$(\tilde{A},\,\tilde\bdb)$ 
~is given in the following corollary of Theorem~\ref{t:ftngs}.

\begin{corollary}
~Let ~$(A,\bdb)$ ~and ~$(\tilde{A},\tilde\bdb)$ ~be as in 
Problem~\ref{p:lsgs} where the underlying linear system 
~$A\,\bdx\,=\,\bdb$ ~is consistent.
~Assume the data matrix ~$\tilde{A}$ ~is sufficiently accurate
with ~$\big\|A-\tilde{A}\big\|_2\,<\,\frac{1}{2}\,\big\|A^\dagger\big\|_2^{-1}$.
~Further assume an error bound ~$\beta\,>\,\big\|A-\tilde{A}\big\|_2$ 
~is attainable
and is sufficiently tight so that 
~$\beta\,<\,\big\|A^\dagger\big\|_2^{-1}-\big\|A-\tilde{A}\big\|_2$.
~Then Problem~\ref{p:lsgs} ~is solvable by calculating 
~$\sol\big(\tilde{A}_\theta,\,\tilde\bdb_\theta\big)$ ~with the error tolerance 
~$\theta\,=\,\beta$ ~where ~$\tilde\bdb_\theta$ ~is the orthogonal projection
of ~$\tilde\bdb$ ~onto ~$\cR(\tilde{A}_\theta)$.
~Furthermore
\begin{align*}
& \dist{\sol(\tilde{A}_\theta,\tilde\bdb_\theta),\, \sol(A,\bdb)} \\
&~~~~\le~  \|A\|_2\,\big\|A^\dagger\big\|_2\cdot
\frac{\sqrt{4\,\|A^\dagger\,\bdb\|_2^2+1}}{
\|A\|_2 -\|A\|_2\,\|A^\dagger\|_2\,\big\|A-\tilde{A}\big\|_2}
\,\big\|(\tilde{A},\,\tilde\bdb)-(A,\,\bdb) \big\|.
\end{align*}
\end{corollary}

We reiterate that the sensitivity of solving ~$A\,\bdx\,=\,\bdb$
~from empirical data ~$(\tilde{A},\tilde\bdb)$ ~is measured by
\[  \|A\|_2\,\|A^\dagger\|_2 ~~=~~ 
\frac{\sg_1(A)}{\sg_r(A)} ~~\approx~~
 \|\tilde{A}_\theta\|_2\,\|\tilde{A}_\theta^\dagger\|_2,
\]
not infinity or ~$\kappa(\tilde{A})$ ~when the underlying matrix ~$A$ ~is 
singular where ~$r\,=\,\rank{A}$.
~Problem~\ref{p:lsgs} may still be difficult if data are inaccurate, if the
intrinsic condition ~$\|A\|_2\,\|A^\dagger\|_2$ ~is large, or if
the window for setting the error tolerance is too narrow.
%
~The general numerical solution can be computed using existing rank-revealing
tools such as \cite{lee-li-zeng,li-zeng-03} 
and UTV/ULV decomposition%
\cite[\S 5.4.6]{golub-vanloan4} in the following template:

\vspace{2mm}
\begin{itemize}
\item[] set the error tolerance ~$\theta$ ~at or slightly above the 
error bound ~$\bt\,\gtrapprox\,\|\Dl A\|_2$
\item[] {\tt if} ~$r\,=\,\ranka{\theta}{A} \approx n$ ~{\tt then}
\begin{itemize}
\item calculate ~$N\,\in\,\C^{n\times (n-r)}$ ~whose columns form an 
orthonormal basis for the numerical kernel ~$\cK(A_\theta)$
\item solve ~$A\,\bdx\,=\,\bdb$ ~for a particular solution ~$\bdx\,=\,\bdx_*$
by any backward accurate method such as
~$\bdx_*\,=\,(A^\h A+\mu^2\,N\,N^\h)^{-1} A^\h\,\bdb$ ~or Tikhonov 
regularization
\item {\tt output} ~$\sol_\theta(A,\,\bdb)\,=\,
\bdx_*+\cR(N)$.
\end{itemize}
\item[] {\tt else}
\begin{itemize}
\item calculate a decomposition 
~$U\,S\,V^\h \,=\, A_\theta$ ~with
~$S\,\in\,\C^{r\times r}$, ~$U^\h U\,=\,I$ ~and ~$V^\h V\,=\,I$ 
\item solve ~$S\,\bdy\,=\,U^\h\,\bdb$ ~for ~$\bdy\,=\,\bdy_*$ ~and
obtain the truncated SVD solution ~$\bdx_*\,=\,V\,\bdy_*$.
\item {\tt output}: ~$\sol_\theta(A,\,\bdb)\,=\,\bdx_*+\cR(V)^\perp$
\end{itemize}
\item[] {\tt end if}
\end{itemize}
\vspace{2mm}

As we shall establish in \S\ref{s:nps},
the particular solution component ~$\bdx_*$ ~of ~$\sol_\theta(A,\bdb)$ ~in
the above template can be computed by any backard accurate numerical algorithm 
including Tikhonov regularization and truncated SVD.
~Computation of general numerical solution is implemented in the 
Matlab package {\sc NAClab} \cite{naclab} as the functionality 
{\tt LinearSolve} (c.f. \cite{solve} and the supplementary material).
~The general guideline for the error tolerance is to set it at or slightly
larger than a known data error bound ~$\bt\,>\,\|A-\tilde{A}\|_2$ ~if the 
application allows such an adjustment.
~We conclude this section with the following example.

\begin{example}\label{e:uf1}\em
~Revisiting the linear system in Example~\ref{e:uf}, the data error bound can 
be estimated as 
~$\|\Dl A\|_2 \,\le\, \|\Dl A\|_{_F}\,\le\,4.5\times 10^{-4}$
~where ~$A$ ~is the underlying matrix since the 
entrywise error bound is ~$5\times 10^{-5}$.
~The error tolerance ~$\theta$ ~can be set at or slightly larger than the 
error bound, say ~$\theta\,=\,0.0005$.
~The numerical solution of the system (\ref{uf123}) within ~$0.0005$ 
~in the affine Grassmannian ~$\cA_7(\C^9)$ ~is a representation of
\begin{align*}
(u_1,\,u_2,\,& u_3) ~~=\\
&(\mbox{\tiny
$~~.90710 + .33322\,x + .71029\,x^2 + .59968\,x^3,~-.79946 + .06694\,x,
~~~\,1.12433 - .06648\,x + .08926\,x^2$}) \\
 + t_1\,&(\mbox{\tiny
$-.27897 - .08391\,x - .17878\,x^2 + .08424\,x^3,~-.35739 - .47261\,x,
~~~~~.12212 - .33612\,x - .63016\,x^2$})\\
+t_2\,&(\mbox{\tiny
$-.21387 + .29319\,x - .18465\,x^2 + .46503\,x^3,~-.55471 + .18011\,x,~
-.46785 + .03542\,x + .24016\,x^2$})
\end{align*}
(c.f. supplementary material).
~The general numerical solution ~$\sol_\theta(\tilde{A},\tilde\bdb)$ 
~is of a healthy sensitivity 
~$\|\tilde{A}_\theta\|_2\,\|\tilde{A}_\theta^\dagger\|_2 \,\approx\, 17.19$,
~not the infinite ~$\kappa(A)$ ~or the large 
~$\kappa(\tilde{A})\,\approx\,2.29\times 10^6$.
~The three components of ~$\sol_\theta(\tilde{A})$
~form an invertible polynomial transformation matrix as shown in (\ref{ufg})
with the numerical inverse
\[  \left[\mbox{\tiny $\begin{array}{l}
.55101 - .91839\,x^2, ~-2.33985 + .70128\,x + .30047\,x^2 + .00001\,x^3,  
~~.71342+ 1.83982\,x + 1.12986\,x^2 + .00002\,x^3 \\
-.36735 + .18367\,x - .73471\,x^5,  
~~-.43135 - 1.09668\,x + .33663\,x^2 - 1.72768\,x^3 + .56105\,x^4
+ 0.24037\,x^5, \\
~~~~~~~~~~~~~~~~~~~~~~~~~~~~~~~~~~~~~~~~~~~~~~~~~~~~~~~~
-.9954 + .44186\,x + .8831\,x^2 + 1.11306\,x^3 + 1.47187\,x^4 + .90389\,x^5 \\
.18366 + .36734\,x^2 + .55103\,x^4, 
~~1.58108 - .53294\,x + 1.17553\,x^2 - .42079\,x^3 - .18027\,x^4, \\
~~~~~~~~~~~~~~~~~~~~~~~~~~~~~~~~~~~~~~~~~~~~~~~~~~~~~~~~~~~~~~~~~~
-1.28338 - 1.39825\,x - 1.28672\,x^2 - 1.10389\,x^3 - .6779\,x^4
\end{array}$}\right].
\]
\end{example}

\begin{remark}\em
~An ~$m\times n$ ~system ~$A\,\bdx\,=\,\bdb$ ~with ~$m\,>\,n$ ~is inconsistent
for almost all ~$\bdb\,\in\,\C^m$ ~and its least squares solution is usually
studied in the literature.
~In fact, the least squares solution can be considered whenever 
~$\bdb\,\not\in\,\cR(A)$ ~even if ~$m\,\le\,n$.
~There are substantial differences between the conventional solution
and the least squares solution.
~In Theorem~\ref{t:ftngs} and throughout this paper, our elaboration is 
restricted to the conventional solution so that 
~$\sol(A,\,\bdb)\,=\,\emptyset$ ~for inconsistent systems and the nonempty
set of least squares solutions is beyond the scope.
~The sensitivity of the least squares solution is well-known 
to be ~$\kappa(A)^2$ (see, e.g.  \cite[\S 20.1]{Higham2}) ~in contrast to 
~$\kappa(A)$ ~for the (conventional) solution in Theorem~\ref{t:ftngs}.
\end{remark}

\section{Particular solution of a singular linear system}\label{s:nps}

There are many applications where only a particular solution is needed among 
the infinitely many solutions of a singular linear system 
~$A\,\bdx\,=\,\bdb$ ~and it makes little difference which particular
solution is obtained. 
~For such applications, the problem of finding  a 
{\em numerical particular solution} can be stated as follows. 

\begin{problem}[Numerical Particular Solution]
\label{p:ls1s}
Assume a linear system ~$A\,\bdx\,=\,\bdb$ ~is consistent
where the entries of ~$A$ ~and ~$\bdb$ ~may be known through empirical
data of limited accuracy. 
~Find a numerical particular solution 
~$\tilde\bdx$ ~that approximates an exact solution 
~$\bdx_*\,\in\,\sol(A,\bdb)$ ~with the error
~$\|\tilde\bdx-\bdx_*\|_2$ ~at an acceptable level.
\end{problem}

There are regularization approaches such as the Tikhonov method 
\cite[\S 6.1.5]{golub-vanloan4}\cite{HansenBook,Neu98}
that can produce 
approximate particular solutions with high backward accuracy.
~For any backward accurate numerical solution ~$\tilde\bdx$ ~of the system 
~$A\,\bdx\,=\,\bdo$ ~in the sense that there is a pair 
~$(\tilde{A},\,\tilde\bdb)$ ~such that ~$\tilde{A}\,\tilde\bdx\,=\,\tilde\bdb$
~and ~$\|(A,\,\bdb)-(\tilde{A},\,\tilde\bdb)\|$ ~is at an acceptable level,
we shall call ~$\tilde\bdx$ ~a {\em numerical particular solution} of 
~$A\,\bdx\,=\,\bdb$.
~The following theorem asserts that every numerical particular solution
approximates one of the exact solutions.

\begin{theorem} \label{t:ls1s}
~Let ~$A\,\in\,\C^{m\times n}$ ~with ~$\rank{A}\,<\,n$ ~and ~$\bdb\,\in\,\cR(A)$.
~Assume ~$\tilde\bdx\,\in\,\C^n$ ~is a backward accurate numerical solution of
~$A\,\bdx\,=\,\bdb$ ~in the sense that
~$\tilde\bdx$ ~is an exact solution of ~$\tilde{A}\,\bdx\,=\tilde\bdb$
~with 
~$\big\|\tilde{A}-A\big\|_2\,\le\,
.46\,\|A^\dagger\|_2^{-1}$.
~Then ~$\tilde\bdx$ ~approximates  an exact solution 
~$\bdx_*\,\in\,\sol(A,\bdb)$ ~with an error bound
\begin{equation} \label{lsee4a}
  \frac{\|\tilde\bdx-\bdx_*\|_2}{\|\bdx_*\|_2}
~~\le~~ 
\frac{\|A\|_2\,\big\|A^\dagger\big\|_2}{1-\big\|A^\dagger\|_2\,\|\Dl A \|_2}\,
\left(
2\,\sqrt{2}\,\frac{\|\Dl A\|_2}{\|A\|_2}
+\frac{\|\Dl\bdb\|_2}{\|\bdb\|_2}\right) 
\end{equation}
assuming ~$\bdb\,\ne\,\bdo$ ~where ~$\Dl A\,=\,A-\tilde{A}$, 
~$\Dl\bdb\,=\,\bdb-\tilde\bdb$, ~or 
\begin{equation} \label{lsee4b}
  \|\tilde\bdx-\bdx_*\|_2 ~~\le~~ 
\frac{\|A\|_2\,\big\|A^\dagger\big\|_2}{1-\big\|A^\dagger\|_2\,\|\Dl A\|_2}\, 
\left(\|\tilde\bdx\|_2\, \frac{\|\Dl A\|_2}{\|A\|_2}
+\frac{\|\Dl\bdb\|_2}{\|A\|_2}\right) 
\end{equation}
if ~$\bdb\,=\,\bdo$.
\end{theorem}

{\em Proof sketch.} ~Let ~$r\,=\,\rank{A}$ ~and ~$\sg_{r+1}(\tilde{A})\,<\,
\theta\,<\,\sg_r(\tilde{A})$. 
~Write ~$\tilde\bdx\,=\,\tilde\bdx_1+\tilde\bdx_2$ ~where 
~$\tilde\bdx_1\,=\,\tilde{A}_\theta^\dagger\,\tilde{A}_\theta\,\tilde\bdx$ ~and
~$\tilde\bdx_2\,=\,(I-\tilde{A}_\theta^\dagger\,\tilde{A}_\theta)\,\tilde\bdx$.
%
~Choose a particular solution 
~$\bdx_*\,=\,A^\dagger\,\bdb + (I-A^\dagger\,A)\,\tilde\bdx_2$ 
~from ~$\sol(A,\bdb)$.
~Since ~$\tilde\bdx_1\,=\,\tilde{A}_\theta^\dagger\,\tilde\bdb$ 
~approximates ~$A^\dagger\,\bdb$, ~$\tilde\bdx_2\,\in\,\cK(\tilde{A}_\theta)$,
~$(I-A^\dagger\,A)\,\tilde\bdx_2\,\in\,\cK(A)$ ~and ~$\cK(\tilde{A}_\theta)$
~approximates ~$\cK(A)$, ~hence ~$\tilde\bdx$ ~is an approximation to the 
particular solution ~$\bdx_*$ ~of ~$A\,\bdx\,=\,\bdb$ ~so the theorem holds.
~Detailed proofs of (\ref{lsee4a}) and (\ref{lsee4b}) are in 
Appendix~\ref{s:prf}. ~\Qed

For the case ~$\bdb\,=\,\bdo$ ~in Theorem~\ref{t:ls1s}, ~the objective is
to solve the homogeneous system ~$A\,\bdx\,=\,\bdo$.
~The inequality (\ref{lsee4b}) includes three cases:

\noindent ~~~Case (i): ~$\tilde\bdb\,=\,\bdo$ ~and ~$\tilde\bdx\,=\,\bdo$. 
~The inequality (\ref{lsee4b}) is trivial and perhaps meaningless 
since ~$\tilde\bdx\,=\,\bdx_*\,=\,\bdo$.
\newline $~~~~$
Case (ii): ~$\tilde\bdb\,=\,\bdo$ ~and ~$\tilde\bdx\,\ne\,\bdo$. 
~Then we can normalize ~$\tilde\bdx$ ~to be a unit vector so that 
(\ref{lsee4b}) becomes
\begin{equation} \label{lsee4d}
\min_{\bdz\in\cK(A)}\,\|\tilde\bdx-\bdz\|_2 
~\le~ \|\tilde\bdx-\bdx_*\|_2 ~\le~
\frac{\|A\|_2\,\big\|A^\dagger\big\|_2}{1-\big\|A^\dagger\|_2\,\|\Dl A\|_2}\, 
\frac{\|\Dl A\|_2}{\|A\|_2}.
\end{equation}
$~~~~$
Case (iii): ~$\tilde\bdb\,\ne\,\bdo$.
~The case is relevant in practical computation by setting the right-hand side
~$\tilde\bdb$ ~as a nonzero random vector of a moderate norm and obtaining a 
numerical particular solution ~$\tilde\bdx$ ~as an exact solution of 
~$\tilde{A}\,\bdx\,=\,\tilde\bdb$ ~with a small ~$\|A-\tilde{A}\|_2$,
~leading to the inverse power iteration.
~The norm ~$\|\tilde\bdx\|_2$ ~is almost always large due to the condition
number ~$\kappa(\tilde{A})\,=\,O(\|A-\tilde{A}\|_2^{-1})$.
~As it turns out pleasantly, the large ~$\|\tilde\bdx\|_2$ ~is exactly 
what is needed as (\ref{lsee4b}) becomes
\begin{equation} \label{lsee4c}
  \left\|
\frac{\tilde\bdx}{\|\tilde\bdx\|_2} - 
\frac{\bdx_*}{\|\tilde\bdx\|_2}\right\|_2
~~\le~~ 
\frac{\|A\|_2\,\big\|A^\dagger\big\|_2}{1-\big\|A^\dagger\|_2\,\|\Dl A\|_2}\,
\frac{1}{\|A\|_2}\, \left(\|\Dl A\|_2 + 
\frac{\|\tilde\bdb\|_2}{\|\tilde\bdx\|_2} \right) 
\end{equation}
The larger the norm ~$\|\tilde\bdx\|_2$ ~achieves, ~the more accurate 
~$\frac{\tilde\bdx}{\|\tilde\bdx\|_2}$ ~is to a particular nontrivial solution of 
the homogeneous system ~$A\,\bdx\,=\,\bdo$.
~Once again, the sensitivity of solving a singular linear system 
~$A\,\bdx\,=\,\bdb$ ~is
~$\|A\|_2\,\big\|A^\dagger\big\|_2 \,=\, \frac{\sg_1(A)}{\sg_r(A)}$,
~not infinity in the sense of finding a numerical particular solution.

Particular solutions of ~$A\,\bdx\,=\,\bdb$ ~can vary arbitrarily but their 
deviations can only stretch in ~$\cK(A)$.
~As the following corollary states, the high sensitivity is near
a direction in ~$\cK(A)$ ~and such a sensitivity may be harmless after all.

\begin{corollary}\label{c:x1x2}
~Let ~$A\,\in\,\C^{m\times n}$ ~with ~$\rank{A}\,<\,n$ ~and ~$\bdb\,\in\,
\cR(A)$.
~Assume ~$\bdx_1$ ~and ~$\bdx_2$ ~are both backward accurate numerical 
particular solutions of ~$A\,\bdx\,=\,\bdb$ ~in the sense that 
~$A_1\,\bdx_1\,=\,\bdb_1$ ~and
~$A_2\,\bdx_2\,=\,\bdb_2$ ~with sufficiently small 
~$\|(A_1,\,\bdb_1)-(A,\,\bdb)\|$ ~and ~$\|(A_2,\,\bdb_2)-(A,\,\bdb)\|$.
~Then there is an ~$\bdx_*\,\in\,\cK(A)$ ~such that 
\begin{align}
\lefteqn{\|(\bdx_1-\bdx_2)-\bdx_*\|_2 ~~\le~~
\|A\|_2\,\big\|A^\dagger\big\|_2\times
} \nonumber \\ &  \times
\left(\frac{\|\bdb-\bdb_1\|_2}{\|A\|_2}+ 
\frac{\|\bdb-\bdb_2\|_2}{\|A\|_2} \right.
 \left. + \frac{\|A-A_1\|_2}{\|A\|_2}\,\|\bdx_1\|_2+ 
\frac{\|A-A_2\|_2}{\|A\|_2}\,\|\bdx_2\|_2 \right).
\label{x1x2}
\end{align}
\end{corollary}

{\em Proof.}  
~Apply the inequality (\ref{lsee4b}) on ~$\tilde{A}\,=\,A$ ~and
~$\tilde\bdx\,=\,\bdx_1-\bdx_2$ ~that satisfies
~$A\,(\bdx_1-\bdx_2) \,=\, (\bdb_1-\bdb)+(\bdb-\bdb_2)+
(A-A_1)\,\bdx_1 + (A_2-A)\,\bdx_2$.
~\Qed

Theorem~\ref{t:ls1s} extends the accuracy result for the inverse iteration in 
spite of the large condition number.
~In \cite{PetWil79}, Peters and Wilkinson described what they called
``exaggerated fears'' in the early days of computer age when the inverse 
iteration 
\begin{equation}\label{invit} (A-\la I)\,\bdx_{k+1} ~~=~~ \bdx_k 
~~~~\mbox{for}~~~ k = 0,\,1,\,\cdots
\end{equation}
at an approximation ~$\la$ ~to an eigenvalue ~$\la_*$ ~of ~$A$ ~was proposed 
for calculating an eigenvector ~$\bdx_*$ ~as a nontrivial 
solution to the homogeneous system ~$(A-\la_* I)\,\bdx\,=\,\bdo$:

\begin{quote}\small
Although [inverse iteration is] basically a simple concept its numerical 
properties have not been widely understood.
~If ~$\la$ ~really is very close to an eigenvalue, the matrix ~$(A-\la I)$
~is almost singular and hence a typical step in the iteration involves the 
solution of a very ill-conditioned set of equations. ...
~The period when inverse iteration was first considered was notable for 
exaggerated fears concerning the instability of direct methods for solving
linear systems and {\em ill-conditioned} systems were a source of particular
anxiety.
~$\cdots$ ~[Few] numerical analysts discuss inverse iteration with any 
confidence.
\end{quote}

It is counterintuitive, and pleasantly surprising nonetheless, 
that ill-condition is not harmful in computing the eigenvector.
~As pointed out in \cite{PetWil79} and by Parlett 
\cite[\S 4.3]{ParlettBook} that errors mainly lie in ~$\cK(A-\la_*\,I)$ ~and are not 
really errors at all:

\begin{quote}\small
[R]oundoff errors can give rise to complely erroneous ``solutions'' to very 
ill-conditioned systems of equations. ...
~Indeed some textbooks have cautioned users not take [$\la$] too close to any 
eigenvalue. ...
Fortunately these fears are groundless and furnish a nice example of confusing
ends with means. ... 
~The error ~$\bde~[\,=\,\bdx_{k+1}-\bdx_*]$, 
~which may be almost as large as the exact solution of
~[$(A-\la I)^{-1}\,\bdx_k$], ~is almost entirely in the direction of ~[the
eigenvector].
...
~The result is alarming if we had hoped for an accurate solution of
[(\ref{invit})] (the means) but is a delight in the search for 
[the eigenvector] (the end).
\end{quote}

Theorem~\ref{t:ls1s} concludes, in fact, that the fears of solving a highly
ill-conditioned linear system may also be exaggerated for non-homogeneous
systems as well when the underlying system ~$A\,\bdx\,=\,\bdb$ 
~is consistent and singular, as long as the numerical solution is 
backward accurate and the intrinsic sensitivity measure 
~$\|A\|_2\,\|A^\dagger\|_2$ ~is moderate.
~The variation between any two numerical particular solutions can be large
but the difference falls harmlessly in the kernel of ~$A$.
~In other words, the ``error'' is actually a part of the solution.

\begin{example}\em
~The system ~$\tilde{A}\,\bdx\,=\,\tilde\bdb$ ~in (\ref{54646}) is a 
representation of the underlying system ~$A\,\bdx\,=\,\bdb$ ~with 
~$\|\Dl A\|_2\,\le\,\|\Dl A\|_{_F}\,\le\,4.5\times 10^{-4}\,=\,\theta$
~from the entrywise error bound ~$0.5\times 10^{-5}$.
~Rounded to five digits after the decimal point,
two numerical particular solutions ~$\tilde\bdx_0$ ~and ~$\bdx_1$
~by turncated SVD ~$\tilde{A}_\theta^\dagger\,\tilde\bdb$ ~and 
Matlab ``$\backslash$'', respectively, are
\begin{eqnarray*} 
\tilde\bdx_0 &~~=~~& [\mbox{\scriptsize
~0.90711,\,  0.33322,\, ~0.71029,\,  0.59968,\, -0.79946,\, ~0.06694,\,  
1.12433,\, -0.06648,\,  ~0.08926}]^\h \\
\tilde\bdx_1 &~~=~~& [\mbox{\scriptsize
-0.78366,\,  0.47296,\, -0.45954,\,  1.83637,\, -3.47453,\, -1.81379,\,  
0.84635,\, -1.57209,\, -2.41843}]^\h
\end{eqnarray*}
with both residuals 
~$\|\tilde{A}\,\tilde\bdx_0 - \tilde\bdb\|\,\approx\, 8.1\times 10^{-5}$
~and
~$\|\tilde{A}\,\tilde\bdx_1 - \tilde\bdb\|\,\approx\, 5.3\times 10^{-5}$
~roughly within the data error bound. 
~The two numerical particular solutions are far apart with 
~$\|\tilde\bdx_0-\tilde\bdx_1\|\,\approx\,5.01$ ~as predicted by the 
large condition number ~$\kappa(\tilde{A})\,\approx\,2.29\times 10^6$. 
~However, the underlying system is 
consistent and singular with a healthy sensitivity 
~$\|A\|_2\,\|A^\dagger\|_2\,\approx\,\|\tilde{A}_\theta\|_2\,
\|\tilde{A}_\theta^\dagger\|_2\,\lessapprox\,17.19$.
~Both ~$\tilde\bdx_0$ ~and ~$\tilde\bdx_1$ ~are accurate approximations to
different exact solutions with estimate error bounds ~$0.00186$ ~and
~$0.00177$ ~respectively, and actual relative errors are 
~$4.49\times 10^{-5}$ ~and ~$0.93\times 10^{-5}$ ~in the same level
of the data error.
\end{example}

\section{Bona fide ill-conditioned linear systems}

A linear system ~$A\,\bdx\,=\,\bdb$ ~is truely ill-conditioned when 
~$\|A\|_2\,\|A^\dagger\|_2$ ~is large regardless of its rank.
~When ~$A$ ~is of full column rank, ~$\bdb\,\in\,\cR(A)$ ~and the condition 
number ~$\kappa(A)$ ~is huge, ~the solution uniquely exists but in general 
can not be computed accurately from perturbed data using whatever algorithm.
~The system is de facto rank-deficient in a practical sense. 
~Even in such cases, a stable general numerical solution may still be 
attainable in an affine Grassmannian from empirical data, and the underlying 
solution can be accurately approximated by a vector in the affine subspace as
the general numerical solution. 

\begin{theorem}\label{t:ics}
~Assume ~$A\,\in\,\C^{m\times n}$, ~$\bdx_*\,\in\,\C^n$ ~and
~$\bdb\,=\,A\,\bdx_*$.
~Let ~$r$ ~be any integer with ~$\sg_r(A)\,>\,\sg_{r+1}(A)$. 
~For any ~$(\tilde{A},\tilde\bdb)\,=\,(A+\Dl A,\, \bdb+\Dl\bdb)$ ~serving as
empirical data of ~$(A,\,\bdb)$ ~with 
\begin{equation}\label{DlA}
 \|\Dl A\|_2~~<~~ \min\{\sg_r(A)-\sg_{r+1}(A),\,(2\,\sqrt{3}-3)\,\sg_r(A)\},
\end{equation}
there is an ~$\tilde\bdx\,\in\,\sol_\theta(\tilde{A},\tilde\bdb)$ 
~with ~$\sg_{r+1}(\tilde{A})\,<\,\theta\,<\,\sg_r(\tilde{A})$ ~such that
\begin{equation}\label{ics}
\frac{\|\tilde\bdx-\bdx_*\|_2}{\|\bdx_*\|_2} ~~\le~~ 
\frac{\sg_1(A)}{\sg_r(A)}\, 
\frac{1}{1-\frac{\sg_{r+1}(A)-\|\Dl A\|_2}{\sg_r(A)}}\,
\left((2+\sqrt{2})\,\frac{\|\Dl A\|_2}{\|A\|_2}+
\frac{\|\Dl\bdb\|_2}{\|\bdb\|_2}\right).
\end{equation}
\end{theorem}

{\em Proof sketch.} 
~Since ~$\bdx_*$ ~is a backward accurate solution of the linear system
~$\tilde{A}_\theta\,\bdx\,=\,\tilde{\bdb}_\theta$, ~Theorem~\ref{t:ls1s}
applies from a reversed perspective. 
~The detailed proof is in Appendix~\ref{s:prf}. ~\Qed

In the following example, the underlying system is 
ill-conditioned but truly nonsingular.  
~All known numerical algorithms including regularization methods 
produce solutions that are inaccuate as single vectors but highly accurate
as the vector component of a general numerical solution that
is perfectly conditioned and contains accurate approximations to 
the underlying exact solution.

\begin{example}\label{e:ics} \em
~Consider the polynomial division problem in the form of the
equation
\[  (x+10)\,q+\rho ~=~ \mbox{\footnotesize $
\frac{1}{3}\,x^8+ 4\,x^7+\frac{23}{3}\,x^6+\frac{34}{3}\,x^5+
15\,x^4+ \frac{56}{3}\,x^3+\frac{67}{3}\,x^2+26\,x+\frac{89}{3}$}
\]
for the quotient ~$q$ ~and ~the constant remainder ~$\rho$.
~There is a unique solution which consists of ~$q \,=\, 
\mbox{\footnotesize $\frac{1}{3}\,(x^7+2\,x^6+\cdots+7\,x+8)$}$ ~and 
~$\rho\,=\,3$.
~The corresponding linear system is of the form ~$A\,\bdx\,=\,\bdb$ ~where
\[  A ~=~ \left[\mbox{\scriptsize $\begin{array}{cccc}
1 &&& \\ 10 & 1 && \\ & \ddots & \ddots & \\ & & 10 & 1
\end{array}$}\right], 
~~~\bdb ~=~ \mbox{$\frac{1}{3}$}\,\left[\mbox{\scriptsize $\begin{array}{c}
1 \\ 12 \\ \vdots \\ 89 \end{array}$}\right]
\]
with the exact solution ~$\bdx_*\,=\,\frac{1}{3}\,[1,2,\cdots,9]^\h$ 
~that is attainable in symbolic computation using the exact data in 
rational number format.
~In Matlab single precision arithmetic, the system is represented as
perturbed data
~$\tilde{A}\,\bdx\,=\,\tilde\bdb$ ~where ~$\tilde{A}\,=\,A$ ~and
~$\tilde\bdb \,=\, [\mbox{\scriptsize 
$.3333333,\,4.0,\,7.6666665,\,11.333333,\,15.0,\,18.666666,\,
22.333334,\,26.0,\,29.666666$}]^\h$.
~The singular values
~{\scriptsize $10.9461079\,>\,10.7891169\,>\,\cdots\,>\, 
9.0683689\,>\,9.9\times 10^{-9}$}
~indicate that it is practically impossible to calculate the single-vector
solution with any meaningful accuracy using such data.
~Table~\ref{t:ics1} shows three sample numerical solutions: ~$\bdx_1$ 
~by a straightforward application of the Matlab command 
{\tt A$\backslash$b}, ~a Tikhonov regularization solution 
~$\bdx_2\,=\,(A^\h\,A+\al^2 I)^{-1}\,A^\h\,\tilde\bdb$ ~at, say ~$\al\,=0.001$
~and the truncated SVD solution 
~$\bdx_3\,=\,A_\theta^\dagger\,\tilde\bdb$ ~with an error tolerance
that is roughly 
~$\theta\,=\,\|\bdb\|_2\,\eps\,\approx\,3.18\times 10^{-6}$ ~where ~$\eps$ ~is
the unit roundoff. 
~As expected from the condition number ~$\kappa(A)\,\approx\,1.1\times 10^9$, 
~none of the solutions can be considered accurate {\em as a single vector}.
~On the other hand, the general numerical solution
~$\sol_\theta(\tilde{A},\tilde\bdb)$ ~is almost perfectly
conditioned at 
~$\big\|\tilde{A}_\theta\big\|_2\,\big\|\tilde{A}_\theta^\dagger\big\|_2
\,\approx\,1.21$. 
~The three solutions ~$\bdx_1$, ~$\bdx_2$ ~and ~$\bdx_3$ ~that are inaccurate 
as individual vectors are all accurate as the component
~$\tilde\bdu$ ~of the general numerical solution 
~$\sol_\theta(\tilde{A},\tilde\bdb)\,=\,\tilde\bdu+\cK(\tilde{A}_\theta)$
~with ~$\cK(\tilde{A}_\theta)\,=\,\spn\{\tilde\bdv\}$ ~where
\begin{equation}\label{v}
 \tilde\bdv ~=~ [\mbox{\tiny $ .0,\, -.0000001,\, .0000010,\,  
-.0000099,\,   .0000995,\,  -.0009950,\,   .0099499,\,  -.0994987,\,.9949875
$}]^\h.
\end{equation}
All ~$\bdx_j+\cK(A_\theta)$ ~for ~$j=1,2,3$ ~are nearly identical in the
affine Grassmannian ~$\cA_1(\C^9)$ ~and each contains a particular vector
~$\hat\bdx_j$ ~that is an accurate approximation to the exact 
solution ~$\bdx_*$ ~as shown in the bottom part of Table~\ref{t:ics1}.
~The errors $\frac{\|\hat\bdx_j-\bdx_*\|_2}{\|\bdx_*\|_2}$ ~are all within
the bound ~$8.28\times 10^{-7}$ ~predicted by (\ref{ics}).

\begin{table}[ht] \scriptsize
\begin{center}
\begin{tabular}{|r|c|c|}\hline
\multicolumn{1}{|c|}{solution} & numerical (single vector) solution & error \\
\multicolumn{1}{|c|}{type} &
with incorrect digits crossed out
& ~$\frac{\|\bdx_j-\bdx_*\|_2}{\|\bdx_*\|_2}$ \\
\cline{1-1}\cline{3-3} 
8 digits of ~~$\bdx_*$  & \tiny \tt .3333333  .6666667  1.0000000  
  1.3333333  1.6666667  2.0000000  2.3333333  2.6666667  3.0000000 
& \\
Matlab ``$\backslash$'' ~~$\bdx_1$ &
 \tiny \tt .3333333  .666666\st{5}  1.00000\st{14}  
  1.3333\st{187}  1.666\st{8129}  1.998\st{5371} 
 2.3\st{479633}  2.\st{5203667}  \st{4.4629993} & \tiny \tt 0.2612916 \\
Tikhonov ~~$\bdx_2$ &  \tiny \tt .3333333  .66666\st{70}  
0.99999\st{71}  1.3333\st{607}  1.666\st{3934} 
 2.00\st{27277}  2.3\st{060572}  2.\st{9394238}  \st{0.2724303} & \tiny \tt 
0.4871437 \\
trunc. SVD ~~$\bdx_3$ &  \tiny \tt .333333\st{5} 0.666666\st{9} 
0.99999\st{67}  1.3333\st{603}  1.666\st{3938} 2.00\st{27270} 
2.3\st{060613} 2.\st{9393935} \st{0.2727296} & \tiny \tt 0.4870902
\\ \cline{2-3}\cline{2-3}
 & particular ~$\hat\bdx_j\,=\,\bdx_j+t_j\,\tilde\bdv\,\in\,
\bdx_j+\cK(A_\theta)$ ~nearest to ~$\bdx_*$
& error \\
 & with ~$\tilde\bdv$ ~in 
(\ref{v}) 
& $\frac{\|\hat\bdx_j-\bdx_*\|_2}{\|\bdx_*\|_2}$
\\ \cline{3-3}
~$\bdx_1+\cK(A_\theta)$ & ~$\hat\bdx_1\,=\,\bdx_1+t_1\,\tilde\bdv
\,\approx\,\bdx_*$ 
~with ~$t_1\,=\,${\tiny \tt ~-1.4703701} & ~$8.8\times 10^{-8}$ \\
~$\bdx_2+\cK(A_\theta)$ & ~$\hat\bdx_2\,=\,\bdx_2+t_2\,\tilde\bdv
\,\approx\,\bdx_*$  
~with ~$t_2\,=~\,${\tiny \tt 2.7413113} & ~$1.5\times 10^{-7}$ \\
~$\bdx_3+\cK(A_\theta)$ & ~$\hat\bdx_3\,=\,\bdx_3+t_3\,\tilde\bdv
\,\approx\,\bdx_*$  
~with ~$t_3\,=~\,${\tiny \tt 2.7410104} & ~$1.9\times 10^{-7}$ \\ \hline
\end{tabular}
\end{center}
\caption{\footnotesize 
For ~$A\,\bdx\,=\,\tilde\bdb$ ~in Example~\ref{e:ics}, 
numerical solutions ~$\bdx_1$, ~$\bdx_2$ ~and ~$\bdx_3$ 
~by Matlab ``$\backslash$'', ~Tikhonov regularization and truncated SVD 
respectively in comparison with the exact solution ~$\bdx_*$, 
~as well as the accuracies of ~$\bdx_1$, ~$\bdx_2$ ~and ~$\bdx_3$ ~as 
a component of the general numerical solution.
}\label{t:ics1}
\end{table}

The linear system in Example~\ref{e:ics} is 
nonsingular in theory but practically underdetermined in numerical computation.
%
%
~Suppose an additional piece of information becomes
available, say the remainder ~$\rho\,=\,3$.
~One can impose such a constraint on the general numerical solution 
~$\{\tilde\bdu+t\,\tilde\bdv ~|~ t\,\in\,\C\}$ ~at the trailing component as
~$0.2727296+.9949875\,t \,=\, 3$, ~obtaining ~$t\,=\,2.7410097$
~corresponding to a numerical solution with a relative error 
~$1.79\times 10^{-7}$ ~in the same order of the data.
\end{example}

%
%

\appendix

\section{Lemmas}\label{s:lem}

\begin{lemma}\label{l:d}
~Let ~$A,\,\tilde{A}\,\in\,\C^{m\times n}$ ~with ~$\Dl A\,=\tilde{A}-A$.
~Assume ~$\sg_r(A)\,>\,\sg_{r+1}(A)$.
\begin{itemize}\parskip-0mm
\item[\em (i)] {\em (Wedin)}
~If ~$\|\Dl A\|_2\,<\,\frac{1}{2}\,\big(\sg_r(A)- \sg_{r+1}(A)\big)$, 
~then
\begin{align} 
\lefteqn{\dist{\cK(A_\theta),\,\cK(\tilde{A_\theta})}} \nonumber \\
&~\le ~\frac{\sg_1(A)}{\sg_r(A)}\,
\frac{1}{1-\frac{\sg_{r+1}(A)+\|\Dl A\|_2}{\sg_r(A)}}\,
	\frac{\big\|\Dl A\big\|_2}{\|A\|_2}
\label{nkdis}  
\\
&~\le~  \frac{\sg_1(A)}{\sg_r(A)}\,
\frac{2}{1-\frac{\sg_{r+1}(A)}{\sg_r(A)}}\,
	\frac{\big\|\Dl A\big\|_2}{\|A\|_2}
\nonumber
\end{align}
for any ~$\theta\,\in\,\big(\sg_{r+1}(A),\,\sg_r(A)\big)\cap
\big(\sg_{r+1}(\tilde{A}),\,\sg_r(\tilde{A})\big)\,\ne\,\emptyset$.
\item[\em (ii)]  ~If ~$\rank{A}\,=\,\rank{\tilde{A}}\,=\,r$ ~and 
~$\|\Dl A\|\,<\,\sg_r(A)$, ~then
\begin{equation} \label{kdis}
\dist{\cK(A),\,\cK(\tilde{A})} ~~\le~~ 
	\frac{\sg_1(A)}{\sg_r(A)}\, \frac{\big\|\Dl A\big\|_2}{\|A\|_2}.
\end{equation}
\end{itemize}
\end{lemma}

{\em Proof}. 
~Assertion (i) is established by Wedin \cite{wedin}
(also see \cite[Theorem 4.4]{StewSun} and \cite[Theorem 3.3]{Han87}).
~To prove (ii), let the singular value decompositions of ~$A$ ~and 
~$\tilde{A}$ ~be
\[ 
A ~=~ [U_1,\,U_2]\,
\left[\begin{array}{cc} \Sigma_1 & \\ & O \end{array}\right]\,[V_1,\,V_2]^\h
~~~~\mbox{and}~~~~
\tilde{A} ~=~ [\tilde{U}_1,\,\tilde{U}_2]\,
\left[\begin{array}{cc} \tilde{\Sigma}_1 & \\ & O\end{array}
\right]\, [\tilde{V}_1,\,\tilde{V}_2]^\h
\]
respectively where ~$\Sigma_1,\,\tilde\Sigma_1\,\in\,\C^{r\times r}$. 
~Then 
\[ -\tilde{V}_2^\h \Dl A^\h ~~=~~ \tilde{V}_2^\h\,\tilde{A}^\h 
-\tilde{V}_2^\h\,\Dl A^\h ~~=~~ \tilde{V}_2^\h\,A^\h \\
~~=~~  (\tilde{V}_2^\h V_1)\,(\Sigma_1^\h\,U_1^\h)
\]
and thus 
\[ \dist{\cK(A),\,\cK(\tilde{A})} ~~=~~
\big\| \tilde{V}_2^\h V_1 \big\|_2 
~~\le~~ \frac{\| \Dl A \|_2}{\sg_r(A)}
~~\le~~ \frac{\sg_1(A)}{\sg_r(A)}\,\frac{\|\Dl A \|_2 }{\|A\|_2}.
~\mbox{\Qed}
\]

\begin{lemma}\label{l:spd}
~Let ~$\cU$ ~be a subspace of ~$\C^n$ ~and ~$U$ ~be a matrix whose columns 
form an orthonormal basis for ~$\cU$.
~For every subspace ~$\cV$ ~of ~$\C^n$ ~of the same dimension as ~$\cU$
~with ~$\dist{\cU,\,\cV}\,<\,1$, ~there is a matrix ~$V$ ~whose columns form 
a basis for ~$\cV$ ~such that ~$\|U-V\|_2 \,\le\, \dist{\cU,\cV}.$ 
\end{lemma}

%
{\em Proof.}
~~Let ~$G$ ~be any matrix whose columns form an 
orthonormal basis for ~$\cV$ ~and ~$[G,\hat{G}]$ ~be a unitary matrix. 
~Then, for any unit vector ~$\bdx\,\in\,\C^n$, 
\[  1 ~~=~~ \big\|[G,\hat{G}]^\h\,U\,\bdx\big\|_2^2 ~~=~~
\|(G^\h\,U)\,\bdx\|_2^2 + \|(\hat{G}^\h\,U)\,\bdx\|_2^2 ~~\le~~
\|(G^\h\,U)\,\bdx\|_2^2 + \dist{\cU,\cV}^2
\]
leading to
~$\|(G^\h\,U)\,\bdx\|_2^2 \,\ge\, 1-\dist{\cU,\cV}^2\,>\,0$,
~implying ~$G^\h\,U$ ~is invertible so that columns of ~$V\,=\,G\,(G^\h\,U)$
~form a basis for ~$\cV$, ~and 
~$\|U-V\|_2 \,=\, \big\|(U\,U^\h- G\,G^\h)\,U\big\|_2$ ~that is 
less than or equals to ~$\dist{\cU,\,\cV}$.  ~\Qed

\begin{lemma}\label{l:A}
~Let ~$A\,\in\,\C^{m\times n}$ ~with ~$\sg_r(A)\,>\,\theta\,>\,\sg_{r+1}(A)$.
%
\begin{itemize}\parskip-0mm
\item[\em (i)] ~For every ~$\mu\,\in\,\big[\sg_r(A),\,\sg_1(A)\big]$, 
~let ~$N\,\in\,\C^{n\times(n-r)}$ ~be a matrix whose columns form an 
orthonormal basis for ~$\cK(A_\theta)$.
~Then
\begin{eqnarray}\label{002}
\mbox{\footnotesize $\left\|\left[\begin{array}{c} \mu\,N^\h \\ A
\end{array}\right]\right\|_2$} & = & 
\max\big\{\sg_1(A),\,\sqrt{\mu^2+\sg_{r+1}(A)^2}\big\} \\ 
& & ~~~~~\in~ \left[ \|A\|_2, ~\sqrt{2}\,\|A\|_2\right) \nonumber \\
\mbox{\footnotesize $\left\|\left[\begin{array}{c} \mu\,N^\h \\ A
\end{array}\right]^\dagger\right\|_2$} & =  &
\max\left\{\frac{1}{\sg_r(A)},\,\frac{1}{\sqrt{\mu^2+\eta^2}}\right\}
~\le~ \|A_\theta^\dagger\|_2
\label{003}
\end{eqnarray}
where ~$\eta\,=\,\sg_n(A)$ ~if ~$m\,\ge\,n$ ~or ~$\eta\,=\,0$ ~otherwise.
\item[\em (ii)]
~Assume columns of ~$N\,\in\,\C^{n\times(n-r)}$ ~span ~$\cK(A_\theta)$.
~For any ~$\mu\,>\,0$, ~let ~$\bdb\,\in\,\C^m$ ~and 
~$\bdx_*$ ~be the
least squares solution of the linear system
\begin{equation}\label{006}
\mbox{\footnotesize $\left[\begin{array}{c} \mu\,N^\h \\ A
\end{array}\right]$}\, \bdx
~~=~~ 
\mbox{\footnotesize $\left[\begin{array}{c} \bdo \\ \bdb \end{array}\right]$}.
\end{equation}
Then ~$\bdx_*\,=\,A^\dagger\bdb$ ~if ~$\rank{A}\,=\,r$ ~or
~$A_\theta\,\bdx_*\,=\,\bdb_\theta$ ~if ~$\rank{A}\,>\,r$ ~where 
~$\bdb_\theta\,=\,A_\theta\,A_\theta^\dagger\,\bdb$
~is the orthogonal projection of ~$\bdb$ ~onto ~$\cR(A_\theta)$. 
Furthermore, 
\begin{equation}\label{004}
\bdx_*-T\,T^\h\,\bdx_* ~~=~~ A_\theta^\dagger\, \bdb_\theta
\end{equation}
for any ~$T\,\in\,\C^{n\times(n-r)}$ ~with 
~$\cR(T)\,=\,\cK(A_\theta)$ ~and ~$T^\h\,T\,=\,I$. 
\end{itemize}
\end{lemma}

{\em Proof.} ~For the case of ~$m\,\ge\,n$, ~we can write ~$A$ ~in its 
singular value expansion
~$A\,=\,\sg_1\,\bdu_1\,\bdv_1^\h+\cdots+\sg_n\,\bdu_n\,\bdv_n^\h$
~and ~$N = [\bdv_{r+1},\cdots,\bdv_n]\,G$ ~where ~$\bdu_1,\,\ldots,\,\bdu_m$
~and ~$\bdv_1,\,\ldots,\,\bdv_n$ ~are left and right singular vectors 
respectively with a unitary matrix ~$G\,\in\,\C^{(n-r)\times(n-r)}$.
~Write ~$\bdx \,=\, x_1\,\bdv_1+\cdots+x_n\,\bdv_n$.
~Then
\begin{eqnarray*}
\lefteqn{\left\|\mbox{\footnotesize $\left[\begin{array}{c} \mu\,N^\h \\ A
\end{array}\right]$}\,\bdx\right\|_2^2 ~~=~~ 
\left\|\mbox{\footnotesize $
\left[\begin{array}{cc} G^\h & \\ & I \end{array}\right]\,
\left[\begin{array}{c} 
\mu\,[ \bdv_{r+1},\cdots,\bdv_n]^\h \\ A
\end{array}\right]$}\,\bdx\right\|_2^2} \\
& = & \sg_1^2\,|x_1|^2+\cdots+\sg_r\,|x_r|^2+(\mu^2+\sg_{r+1}^2)
\,|x_{r+1}|^2+\cdots+ (\mu^2+\sg_n^2)\,|x_n|^2
\end{eqnarray*}
whose extrema subject to
~$\|\bdx\|_2\,=\,1$ ~are
~$\max\{\sg_1^2,\,\mu^2+\sg_{r+1}^2\}$ ~and
~$\min\{\sg_r^2,\,\mu^2+\sg_n^2\}$, 
leading to (\ref{002}) and (\ref{003}) in the assertion (i).
~The case ~$m\,<\,n$ ~is similar. 

To prove (ii), write the singular value decomposition 
~$A\,=\,U_1\,\Sigma_1\,V_1^\h+U_2\,\Sigma_2\,V_2^\h$ ~where 
~$\Sigma_1\,\in\,\C^{r\times r}$ ~and 
~$\Sigma_2\,\in\,\C^{(m-r)\times (n-r)}$.
~Then ~$\bdx_*$ ~is the solution of the normal equation
~$\mu^2 N N^\h \bdx_* + A^\h A \bdx_* - A^\h \bdb \,=\, \bdo$.
~Namely, we have an orthogonal decomposition
\begin{align}\label{007}
& (V_1\, \Sigma_1^\h\, \Sigma_1\, V_1^\h\, \bdx_* 
- V_1\, \Sigma_1^\h\, U_1^\h\, \bdb)+ 
(V_2\, \Sigma_2^\h\, \Sigma_2\, V_2^\h\, \bdx_* - 
V_2\, \Sigma_2^\h\, U_2^\h\, \bdb +\mu^2\, N\, N^\h\, \bdx_*)  \\
& ~~=~~ \bdo,  \nonumber
\end{align}
implying ~$V_1\, \Sigma_1^\h\, \Sigma_1\, V_1^\h\, \bdx_* 
- V_1\, \Sigma_1^\h\, U_1^\h\, \bdb \,=\, \bdo$ ~and thus 
~$V_1^\h\,\bdx_*\,=\,\Sigma_1^{-1}\,U_1^\h\,\bdb$.
~Since ~$\bdx_* \,=\, V_1\, V_1^\h\, \bdx_* + V_2\, V_2^\h\, \bdx_*$, ~hence
~$A_\theta\,\bdx_* \,=\, (U_1\,\Sigma_1\,V_1^\h)\,(V_1\,V_1^\h\,\bdx_*)
\,=\, U_1\,U_1^\h\,\bdb \,=\, \bdb_\theta$.
Namely ~$\bdx_*$ ~is a particular solutions of the 
system ~$A_\theta\,\bdx = \bdb_\theta$.
~Also,
\[ A_\theta^\dagger\,\bdb_\theta ~=~ 
(V_1\,\Sigma_1^{-1}\,U_1^\h)\,(U_1\,U_1^\h\,\bdb) ~=~
(V_1\, \Sigma_1^{-1}\, U_1^\h)\,\bdb ~=~ V_1\, V_1^\h \bdx_* 
~=~ (I-T\,T^\h)\, \bdx_*.
\]
Finally, if ~$\rank{A}\,=\,r$, ~then ~$A_\theta\,=\,A$ ~so ~$\Sigma_2\,=\,O$
~in (\ref{007}), implying ~$N^\h\,\bdx_*\,=\,\bdo$.
~Consequently ~$\bdx_*\,=\,A^\dagger\,\bdb$ ~from ~(\ref{004}). ~\Qed

The following lemma is a variation of 
Theorem 5.1 in \cite{wedin73} by Wedin and its extension in 
Theorem~3.4 in \cite{Han87} by Hansen.

\begin{lemma}[Wedin, Hansen]\label{l:mns}
~Let ~$A\,\in\,\C^{m\times n}$ ~and ~$\bdb\,\in\,\C^m$.
~Assume, ~for a ~$\theta\,>\,0$, ~we have
~$\ranka{\theta}{A}\,=\,r$ ~and ~$\|\bdb-A_\theta\,A_\theta^\dagger\,\bdb\|_2
\,<\,\theta$. 
~There is a constant
\begin{equation}\label{zeta}
\zeta ~~=~~ 
\|A_\theta^\dagger\,\bdb\|_2 + 
\frac{1+\|A_\theta^\dagger\,\bdb\|_2}{1-\|A_\theta^\dagger\|_2\,
\|A-A_\theta\|_2}
\end{equation}
such that, ~for any ~$\tilde{A}\,=\,A+\Dl A\,\in\,\C^{m\times n}$ ~and 
~$\tilde\bdb\,=\,\bdb+\Dl\bdb\,\in\,\C^m$ ~with
\begin{equation}\label{aad} \|\Dl A\|_2 ~~<~~ \min\left\{ 
\mbox{$\frac{1}{2}$}\,\big(\sg_r(A)-\sg_{r+1}(A)\big),~\sg_r(A)-\theta,
~\theta-\sg_{r+1}(A) \right\},
\end{equation}
the following inequality holds:
\begin{equation}
\big\|A_\theta^\dagger \bdb - \tilde{A}_\theta^\dagger \tilde\bdb\big\|_2
~~\le~~
\frac{\sg_1(A)}{\sg_r(A)}
\left(\zeta\,\frac{\|\Dl A\|_2}{\|A\|_2}\,
+ \frac{\|\Dl\bdb\|_2}{\|A\|_2} \right) + O(\|(\Dl A,\,\Dl\bdb)\|^2).
\label{mnscont}
\end{equation}
As a special case, 
further assume ~$\rank{A}\,=\,r$ ~and ~$\bdb\,\in\,\cR(A)$.
~Then
\begin{equation}
\big\|A^\dagger \bdb - \tilde{A}_\theta^\dagger \tilde\bdb\big\|_2
~~\le~~
\frac{\sg_1(A)}{\sg_r(A)}\,\frac{1}{1-\frac{\|\Dl A\|_2}{\sg_r(A)}}\,
\left(2\,\|A^\dagger\,\bdb\|_2\,\frac{\|\Dl A\|_2}{\|A\|_2}\,
+ \frac{\|\Dl\bdb\|_2}{\|A\|_2} \right).
\label{mnscont1}
\end{equation}
\end{lemma}

{\em Proof.}
The assumption ~$\ranka{\theta}{A}\,=\,r$ ~implies 
~$\sg_{r+1}(A)\,<\,\theta\,<\,\sg_r(A)$ ~and thus
~$\sg_{r+1}(\tilde{A})\,<\,\theta\,<\,\sg_r(\tilde{A})$ ~following (\ref{aad})
so that ~$\ranka{\theta}{\tilde{A}}\,=\,r$ ~as well.
~Then it is straightforward to verify (\ref{mnscont}) from 
the inequality (26a) in \cite{Han87} using 
\[ \|A\,(A_\theta^\dagger\,\bdb)-\bdb\|_2 ~~=~~ 
\|A_\theta\,A_\theta^\dagger\,\bdb-\bdb\|_2 
~~<~~ \theta ~~<~~ \sg_r(A).
\]
The inequality (\ref{mnscont1}) follows from 
\cite[inequality (27a)]{Han87} 
and ~$\bdb\,\in\,\cR(A)$.
~\Qed

\begin{lemma}\label{l:sd}
~At any ~$(A,\,\bdb)\,\in\,\C^{m\times n}\times\C^n$ ~and 
~$\theta\,>\,0$ ~within which ~$\sol_\theta(A,\bdb)$ 
~is well-defined, there is a ~$\dl\,>\,0$ ~such that
~$\sol_\theta(A+\Dl A,\bdb+\Dl\bdb)$ ~is well-defined with the same dimension 
as ~$\sol_\theta(A,\bdb)$ ~if ~$\|(\Dl A,\,\Dl\bdb)\|\,<\,\dl$. 
\end{lemma}

{\em Proof.} 
~Write ~$\tilde{A}\,=\,A+\Dl A$ ~and ~$\tilde\bdb\,=\,\bdb+\Dl\bdb$.
~Since ~$r\,=\,\ranka{\theta}{A}$ ~is well-defined, we have
~$\sg_{r+1}\,<\,\theta\,<\,\sg_r(A)$. 
~Thus ~$\|\Dl A\|_2\,<\,\min\{\sg_r-\theta,\,\theta-\sg_{r+1}\}$ 
~ensures ~$\ranka{\theta}{\tilde{A}}\,=\,r$.
~Let ~$P\,=\,I-A_\theta\,A_\theta^\dagger$ ~and 
~$\tilde{P}\,=\,I-\tilde{A}_\theta\,\tilde{A}_\theta^\dagger$. 
%
~Then 
\[ \tilde{A}-\tilde{A}_\theta ~~=~~ \tilde{P}\,\tilde{A} ~~=~~
\tilde{P}\,\Dl A +(\tilde{P}-P)\,A+ (A-A_\theta)
\]
and by (\ref{nkdis}), 
\[ \|P-\tilde{P}\|_2 ~~=~~
\dist{\cR(A_\theta),\,\cR(\tilde{A}_\theta)} ~~\le~~
\eta\,\frac{\|\Dl A\|_2}{\|A\|_2}
\]
where, assuming 
~$\|\Dl A\|_2\,\le\,\frac{1}{2}\,(\sg_r(A)-\sg_{r+1}(A))$,
\[ \eta ~~=~~ 
\frac{\sg_1(A)}{\sg_r(A)}\,\frac{2}{1-\frac{\sg_{r+1}(A)}{\sg_r(A)}}
~~=~~ 
\frac{2\,\|A_\theta\|_2\,\|A_\theta^\dagger\|_2}{
1-\|A_\theta^\dagger\|_2\,\|A-A_\theta\|_2}
\]
implying
\[ \|A-A_\theta\|_2 - (\eta+1)\,\|\Dl A\|_2 ~~\le~~ 
\|\tilde{A}-\tilde{A}_\theta\|_2 ~~\le~~
\|A-A_\theta\|_2 + (\eta+1)\,\|\Dl A\|_2
\]
and similarly
\begin{align*} \lefteqn{
\|\bdb-\bdb_\theta\|_2 - \eta\,\frac{\|\Dl A\|_2}{\|A\|_2}\,\|\bdb\|_2 
-\|\Dl\bdb\|_2 ~~\le~~ 
\|\tilde{\bdb}-\tilde\bdb_\theta\|_2~~~~~~~~~~~~~~~}\\
& ~~~~~~~~~~~~~~\le~~
\|\bdb-\bdb_\theta\|_2 + \eta\,\frac{\|\Dl A\|_2}{\|A\|_2}\,\|\bdb\|_2 
+\|\Dl\bdb\|_2
\end{align*}
where ~$\bdb_\theta\,=\,\bdb-P\,\bdb$ ~and ~$\tilde\bdb_\theta\,=\,\tilde\bdb
-\tilde{P}\,\tilde\bdb$. 
%
~If ~$\sol_\theta(A,\bdb)$ ~is empty, then 
~$\|(A,\bdb)-(A_\theta,\bdb_\theta)\|\,>\,\theta$ ~and ~thus
~$\|(\tilde{A},\tilde\bdb)-(\tilde{A}_\theta,\tilde\bdb_\theta)\|\,>\,\theta$ 
~when ~$\|(\Dl A,\,\Dl\bdb)\|$ ~is sufficiently small so that 
~$\sol_\theta(\tilde{A},\tilde\bdb)\,=\,\emptyset$ ~as well.
~When ~$\|(A,\bdb)-(A_\theta,\bdb_\theta)\|\,<\,\theta$ ~and
~$\|(\Dl A,\,\Dl\bdb)\|$ ~is sufficiently small, ~we also have
~$\|(\tilde{A},\tilde\bdb)-(\tilde{A}_\theta,\tilde\bdb_\theta)\|\,<\,\theta$ 
~and ~$\sg_{r+1}(\tilde{A})\,<\,\theta\,<\,\sg_r(\tilde{A})$
~so that ~$\sol_\theta(\tilde{A},\tilde\bdb)\,=\,
\tilde{A}_\theta^\dagger\,\tilde\bdb_\theta+\cK(\tilde{A}_\theta)$ 
~has the identical dimension ~$n-r$ ~as ~$\sol_\theta(A,\bdb)$. ~\Qed

\section{Proofs of theorems and corollaries}\label{s:prf}

{\em Proof of Theorem~\ref{t:lsps}.} (p.\,\pageref{t:lsps})
~Let ~$N\,\in\,\C^{n\times (n-r)}$ ~whose columns form an orthonormal basis for
~$\cK(A)$.
~By Lemma~\ref{l:spd}, ~there is an ~$\tilde{N}\,\in\,\C^{n\times (n-r)}$ 
~whose columns form a basis for ~$\cK(\tilde{A})$ ~such that 
~$\|N-\tilde{N}\|_2\,\le\,\dist{\cK(A),\,\cK(\tilde{A})}$.
~For ~$\zeta = \sg_r(A)$, ~denote 
~$B \,=\, \mbox{\scriptsize 
$\left[ \begin{array}{c} \zeta N^\h \\ A \end{array} \right]$}$
~and ~$\tilde{B}\,=\,\mbox{\scriptsize 
$\left[ \begin{array}{c} \zeta \tilde{N}^\h \\ \tilde{A}\end{array} \right]$}$.
~Then ~$A^\dagger\,\bdb\,=\,B^\dagger\,\mbox{\scriptsize 
$\left[ \begin{array}{c} \bdo \\ \bdb \end{array} \right]$}$ ~and
~$\tilde{A}^\dagger\,\tilde\bdb\,=\,\tilde{B}^\dagger\,\mbox{\scriptsize 
$\left[ \begin{array}{c} \bdo \\ \tilde\bdb \end{array} \right]$}$ 
~by Lemma~\ref{l:A} part (ii).
~By ~$\|B-\tilde{B}\|_2\,\le\,\sqrt{2}\,\|\Dl A\|_2$ ~from Lemma~\ref{l:d} 
part (ii), \cite[Theorem 1.4.6, page 30]{bjorck96} and Lemma~\ref{l:A} part (i),
\begin{align}
\left\|B^\dagger\,\mbox{\scriptsize 
$\left[ \begin{array}{c} \bdo \\ \bdb \end{array} \right]$}- 
\tilde{B}^\dagger\,\mbox{\scriptsize 
$\left[ \begin{array}{c} \bdo \\ \tilde\bdb \end{array} \right]$}\right\|_2
&~\le~ 
\frac{\sg_1(B)}{\sg_n(B)}\,\frac{1}{1-\frac{\|B-\tilde{B}\|_2}{\sg_n(B)}}\,
\left(\|\bdx_*\|_2\frac{\|B-\tilde{B}\|_2}{\|B\|_2}+\frac{\|\Dl \bdb\|_2}{\|B\|_2}
\right) \nonumber \\
&~\le~ 
\frac{\sg_1(A)}{\sg_r(A)}\,\frac{1}{1-\frac{\sqrt{2}\,\|\Dl A\|_2}{\sg_r(A)}}\,
\frac{\sqrt{2}\,\|\bdx_*\|_2\,\|\Dl A\|_2+\|\Dl \bdb\|_2}{\|A\|_2}
\label{26428}
\end{align}
leading to (\ref{afsest}).  ~\Qed

{\em Proof Theorem~\ref{t:ftngs}.} (p.\,\pageref{t:ftngs})
~The assertion (i) is true because ~$A_\theta\,=\,A$ 
~and ~$\bdb_\theta\,=\,\bdb$ ~for ~$\theta\,\in\, (0,\sg_r(A))$.
~If ~$\bdb\,\in\,\cR(A)$, ~then 
\[ \sol_\theta(A,\bdb) ~=~ \sol(A,\bdb) ~=~ A^\dagger\bdb+\cK(A). 
\]
Otherwise 
~$\sol_\theta(A,\bdb)\,=\,\sol(A,\bdb) \,=\,\emptyset$ ~if
~$\theta\,<\,\min\big\{\sg_r(A),\,\|A\,A^\dagger\,\bdb-\bdb\|_2\big\}$. 
~The assertion (ii) directly follows from Lemma~\ref{l:mns} and 
Lemma~\ref{l:sd} with 
\[ \xi ~~=~~ \|A_\theta\|_2\,\|A_\theta^\dagger\|_2\,
\frac{\sqrt{\zeta^2+1}\,
}{
\|A\|_2- \|A\|_2\,\|A_\theta^\dagger\|_2\,\|A-A_\theta\|_2}\,+\eps
\]
for any ~$\eps\,>\,0$.

We now prove the assertion (iii), part (a).
~Let ~$\tilde\bdb_\theta\,=\,
\tilde{A}_\theta\,\tilde{A}_\theta^\dagger\,\tilde\bdb$, 
~$P\,=\,I-A\,A^\dagger$ ~and 
~$\tilde{P}\,=\,I-\tilde{A}_\theta\,\tilde{A}_\theta^\dagger$.
~From ~$\|\tilde{A}-\tilde{A}_\theta\|_2$ \,=\, 
$\min_{\rank{B}=r}\|\tilde{A}-B\|_2 \,\le\, \|\Dl A\|_2$,
we have
\begin{align*} \|\tilde\bdb-\tilde\bdb_\theta\|_2 
& ~=~  \|\tilde{P}\,\tilde\bdb\|_2 ~=~  \|\tilde{P}\,\tilde\bdb-P\,\bdb\|_2 
~\le~  
\|\tilde{P}\|_2\,\|\tilde\bdb-\bdb\|_2+\|\tilde{P}-P\|_2\,\|\bdb\|_2 \\
&~\le~ 
\|\Dl \bdb\|_2+ \dist{\cR(\tilde{A}_\theta),\,\cR(A)}\,\|\bdb\|_2 \\
&~\le~  \|\Dl \bdb\|_2+ \frac{\sg_1(A)}{\sg_r(A)}
\frac{2\,\|\bdb\|_2}{\|A\|_2}\, \|\Dl A\|_2  \tag{by (\ref{nkdis})}\\
&~\le~  \sqrt{4\,\|A^\dagger\|_2^2\,\|\bdb\|_2^2+1}\,
\left\|(\Dl A,\,\Dl\bdb) \right\| 
~=~ \sqrt{\omega^2-1}\,\left\|(\Dl A,\,\Dl\bdb) \right\|
\end{align*}
and
\begin{align}\label{654} 
\big\|(\tilde{A},\,\tilde\bdb)-(\tilde{A}_\theta,\,\tilde\bdb_\theta)\big\|
&~~\le~~ \sqrt{\|\Dl A\|_2^2+(\omega^2-1)\,(\|\Dl A\|_2^2+\|\Dl\bdb\|_2^2)}
\nonumber \\
&~~\le~~ \omega\,\big\|(\Dl A,\,\Dl\bdb)\big\| 
~~<~~ \sg_r(A)-\big\|(\Dl A,\,\Dl\bdb)\big\| 
\end{align}
Then, for any ~$\theta$ ~satisfying (\ref{tintv}), 
\[ \sg_{r+1}(\tilde{A}) ~\le~ \big\|(\Dl A,\,\Dl\bdb)\big\|
~<~ \theta ~<~\sg_r(A)-\big\|(\Dl A,\,\Dl\bdb)\big\|
~\le~ \sg_r(\tilde{A})
\]
so ~$\ranka{\theta}{\tilde{A}}\,=\,r$, 
~$\big\|(\tilde{A},\,\tilde\bdb)-(\tilde{A}_\theta,\,\tilde\bdb_\theta)\big\|
\,<\,\theta$ ~and thus ~$\sol_\theta(\tilde{A},\,\tilde\bdb)$ ~is of the same
dimension as ~$\sol(A,\,\bdb)$.
~Since ~$\sol_\theta(\tilde{A},\,\tilde\bdb)\,=\,
\sol(\tilde{A}_\theta,\tilde\bdb_\theta)$, 
~the backward error of ~$\sol_\theta(\tilde{A},\,\tilde\bdb)$ 
~is bounded above by ~$\omega\,\big\|(\Dl A,\,\Dl\bdb)\big\|$ 
~from (\ref{654}).
~Thus (\ref{ngserr}) follows from (\ref{nkdis}) in Lemma~\ref{l:d} and 
(\ref{mnscont1}) in Lemma~\ref{l:mns}, 
leading to the assertion (iii).

We now prove the assertion (b) of part (iii).
~If ~$\sol(A,\bdb)$ ~is empty, ~then ~$\sol_\theta(A,\bdb)\,=\,\emptyset$
~whenever ~$\theta\,<\,\min\{\sg_r(A),\,\|\bdb-A\,A^\dagger\,\bdb\|_2\}$.
~By Lemma~\ref{l:sd}, there is a ~$\dl_1\,>\,0$ ~such that 
~$\sol_\theta(\tilde{A},\,\tilde\bdb)\,=\,\emptyset$ ~for every
~$(\tilde{A},\,\tilde\bdb)$ ~with
~$\big\|(\tilde{A},\,\tilde\bdb)-(A,\,\bdb)\big\|\,<\,\dl_1$.
~Thus ~$\sol_\theta(\tilde{A},\,\tilde\bdb)\,=\, \sol(A,\bdb)$ 
~with both backward and forward errors as zero.  ~\Qed

{\em Proof of Theorem~\ref{t:ls1s}.} (p.\,\pageref{t:ls1s})
~Let 
~$\tilde{A} \,=\, \tilde{U}_1\,\tilde{\Sigma}_1\,\tilde{V}_1^\h+
\tilde{U}_2\,\tilde{\Sigma}_2\,\tilde{V}_2^\h$ ~be the singular value 
decomposition where ~$\tilde\Sigma_1$ ~is ~$r\times r$ ~with ~$r\,=\,\rank{A}$.
~Then ~$\tilde\bdx$ ~is a solution to ~$\tilde{A}\,\bdx\,=\,\tilde\bdb$ 
~implies ~$\tilde{U}_1\, \tilde\Sigma_1\,\tilde{V}_1^\h\, 
\tilde\bdx_1 \,=\, \tilde{U}_1\,\tilde{U}_1^\h\,\tilde\bdb$
~and ~$\tilde{U}_2\, \tilde\Sigma_2\,\tilde{V}_2^\h\, \tilde\bdx_2 
\,=\, \tilde{U}_2\,\tilde{U}_2^\h\,\tilde\bdb$
~where ~$\tilde\bdx_1\,=\,\tilde{V}_1\,\tilde{V}_1^\h\,\tilde\bdx$
~and ~$\tilde\bdx_2\,=\,\tilde{V}_2\,\tilde{V}_2^\h\,\tilde\bdx$. 
~Then ~$\tilde\bdx\,=\,\tilde\bdx_1+\tilde\bdx_2$ ~with
~$\tilde\bdx_1\,=\,\tilde{A}_\theta^\dagger\,\tilde\bdb$
~for any ~$\theta$ ~between
~$\sg_{r+1}(\tilde{A})$ ~and ~$\sg_r(A)-\|\Dl A\|_2$.
~By Lemma~\ref{l:mns} with ~$\hat\bdx\,=\,A^\dagger\, \bdb$, ~we have
\[
\big\|\tilde\bdx_1 - \hat\bdx\big\|_2
 ~~\le~~ 
\frac{\sg_1(A)}{\sg_r(A)}\,
\frac{\|\hat\bdx\|_2}{1-\frac{\|\Dl A\|_2}{\sg_r(A)}}\, 
\left(2\,
\frac{\|\Dl A\|}{\|A\|_2}+ \frac{\|\Dl\bdb\|_2}{\|\bdb\|_2}
\right)
\]
Let columns of ~$N$
~form an orthonormal basis for ~$\cK(A)$.
~Since ~$\tilde\bdx_2\,\in\,\cK(\tilde{A}_\theta)$, 
\begin{eqnarray} 
\|N\,N^\h\,\tilde\bdx_2-\tilde\bdx_2\|_2 
& ~=~  &
\min_{\bdu\in\cK(A)}\|\bdu-\tilde\bdx_2\|_2 \nonumber \\ 
&~~\le~~& 
\dist{\cK(A_\theta),\,\cK(A)}\, \|\tilde\bdx_2\|_2 \nonumber \\
& \le & 
\frac{\sg_1(A)}{\sg_r(A)}\,
\frac{1}{1-\frac{\|\Dl A\|_2}{\sg_r(A)}}\, 
\frac{\|\Dl A\|_2}{\|A\|_2}\,
\|\tilde\bdx_2\|_2  \label{1255}
\end{eqnarray}
by Lemma~\ref{l:d} which, combined with ~$\|\Dl A\|_2\,\le\,
.46\,\sg_r(A) \,<\,(2\,\sqrt{3}-3)\,\sg_r(A) 
$, ~implies 
$\dist{\cK(A_\theta),\,\cK(A)}\,<\,\frac{\sqrt{3}}{2}$ ~and thus
\begin{eqnarray*}  
\|N\,N^\h\,\tilde\bdx_2\|_2 &~~=~~& \|N^\h\,\tilde\bdx_2\|_2 ~~=~~ 
\|(N^\h\,\tilde{V}_2)\,\tilde{V}_2^\h\,\tilde\bdx_2\|_2 \\
& \ge & \sqrt{1-\dist{\cK(A_\theta),\,\cK(A)}^2}\,
\|\tilde{V}_2^\h\,\tilde\bdx_2\|_2 
~~\ge~~ \mbox{$\frac{1}{2}$}\,\|\tilde\bdx_2 \|_2.
\end{eqnarray*}
Let ~$\bdx_*\,=\,\hat\bdx+N\,N^\h\,\tilde\bdx_2$.
~Then ~$\bdx_*$ ~is a particular solution to ~$A\,\bdx\,=\,\bdb$ ~and
~$\|\bdx_*\|_2^2\,=\,\|\hat\bdx\|_2^2+\|N\,N^\h\,\tilde\bdx_2\|_2^2$.
~We have
\begin{eqnarray*}
\lefteqn{\|\tilde\bdx-\bdx_*\| ~~\le~~ \|\tilde\bdx_1-\hat\bdx\|_2+
\|\tilde\bdx_2-N\,N^\h\,\tilde\bdx_2\|_2}\\
& ~\le~ & 
\frac{\sg_1(A)}{\sg_r(A)}\,
\frac{1}{1-\frac{\|\Dl A\|_2}{\sg_r(A)}}\, 
\left( \frac{\|\Dl A\|}{\|A\|_2}
\big(2\,\|\hat\bdx\|_2 + 2 \|N\,N^\h\,\tilde\bdx\|_2\big)
+\|\bdx_*\|\,\frac{\|\Dl \bdb\|_2}{\|\bdb\|_2} \right) \\
& \le & 
\frac{\sg_1(A)}{\sg_r(A)}\,
\frac{\|\bdx_*\|_2}{1-\frac{\|\Dl A\|_2}{\sg_r(A)}}\, 
\left(2\,\sqrt{2}\,\frac{\|\Dl A\|_2}{\|A\|_2}
+\frac{\|\Dl \bdb\|_2}{\|\bdb\|_2}\right) 
\end{eqnarray*}
leading to (\ref{lsee4a}). 
~For the case ~$\bdb\,=\,\bdo$, ~the bound (\ref{lsee4b}) follows from
(\ref{1255})
\[
\|\tilde\bdx_1\|_2 ~~\le~~ 
\frac{\|U_1^\h\,\tilde\bdb\|_2}{\sg_r(\tilde{A})} 
~~\le~~ \frac{\sg_1(A)}{\sg_r(A)}\,\frac{1}{1-\frac{\|\Dl A\|_2}{\sg_r(A)}}\, 
\frac{\|\Dl \tilde\bdb\|_2}{\|A\|_2}
~~~~~\mbox{\Qed}
\]

{\em Proof of Theorem~\ref{t:ics}.}  (p.\,\pageref{t:ics})
~Let ~$\,U_1\,\Sigma_1\,V_1^\h+U_2\,\Sigma_2\,V_2^\h$
~and ~$\tilde{U}_1\,\tilde{\Sigma}_1\,\tilde{V}_1^\h+
\tilde{U}_2\,\tilde{\Sigma}_2\,\tilde{V}_2^\h$ 
~be singular value decompositions of ~$A$ ~and ~$\tilde{A}$ ~respectively
where ~$\Sigma_1,\,\tilde\Sigma_1\,\in\,\C^{r\times r}$.
~Denote ~$A_1\,=\,U_1\,\Sigma_1\,V_1^\h$, 
~$\tilde{A}_1\,=\,\tilde{U}_1\,\tilde{\Sigma}_1\,\tilde{V}_1^\h$,
~$\bdx_1\,=\,A_1^\dagger\,\bdb$, ~$\bdx_2\,=\,\bdx_*-\bdx_1$,
~$\tilde\bdx_1\,=\,\tilde{A}_1^\dagger\,\tilde\bdb$ ~and 
~$\bdr\,=\, A\,\bdx_1-\bdb$. 
~Then, with ~$\bdr\,=\, U_2\,U_2^\h\,\bdb \,=\, U_2\,\Sigma_2\,V_2^\h\,\bdx_2$,
\begin{align*}
\tilde\bdx_1-\bdx_1 &~~=~~ \tilde{A}_1^\dagger\,(\bdb+\Dl\bdb)-\bdx_1 
~~=~~ \tilde{A}_1^\dagger\,(A\,\bdx_1-\bdr+\Dl\bdb)-\bdx_1 \\
&~~=~~ \tilde{A}_1^\dagger\,((\tilde{A}-\Dl A)\,\bdx_1-\bdr+\Dl\bdb)-\bdx_1 \\
&~~=~~ \tilde{A}_1^\dagger\,(-\Dl A\,\bdx_1-\bdr+\Dl\bdb)-
(I-\tilde{A}_1^\dagger\,\tilde{A}_1)\,\bdx_1 \\
&~~=~~ \tilde{A}_1^\dagger\,(-\Dl A\,\bdx_1+\Dl\bdb)
-\tilde{V}_1\,\tilde\Sigma_1^{-1}\,\tilde{U}_1^\h\,
U_2\,\Sigma_2\,V_2^\h\,\bdx_2
-\tilde{V}_2\,\tilde{V}_2^\h\,\bdx_1,
\end{align*}
leading to
\begin{align*}
\lefteqn{
\|(\tilde\bdx_1 +\tilde{V}_2\,\tilde{V}_2^\h\,\bdx_1)-\bdx_1\|_2}\\
 &~~\le~~ \|\tilde{A}_1^\dagger\|_2\,(\|\Dl A\|_2\,\|\bdx_1\|_2 
+\|\Dl\bdb\|_2)
+ \|\tilde\Sigma_1^{-1}\|_2\,\|\Sigma_2\|_2\,\|\tilde{U}_1^\h\,U_2\|_2\,
\|\bdx_2\|_2 \\
 &~~\le~~ \|\tilde{A}_1^\dagger\|_2\,(\|\Dl A\|_2\,\|\bdx_1\|_2 
+\|\Dl\bdb\|_2)
+ \dist{\cR(U_1),\,\cR(\tilde{U}_1)}\,\|\bdx_2\|_2 \\
&~~\le~~ 
\frac{\sg_1(A)}{\sg_r(A)}\, 
\frac{\|\bdx_*\|_2}{1-\frac{\sg_{r+1}+\|\Dl A\|_2}{\sg_r(A)}}\,
\left(\sqrt{2}\,\frac{\|\Dl A\|_2}{\|A\|_2}+\frac{\|\Dl\bdb\|_2}{\|\bdb\|_2}
\right).  
\end{align*}
Let ~$\tilde\bdx\,=\,\tilde\bdx_1+\tilde\bdx_2\,\in\,
\sol_\theta(\tilde{A},\,\tilde\bdb)$ ~with
\[ \tilde\bdx_2 ~~=~~ \tilde{V}_2\,\tilde{V}_2^\h\,\bdx_1
+\tilde{V}_2\,\tilde{V}_2^\h\,\bdx_2\,\in\,\cK(\tilde{A}_\theta).
\]
Then
\begin{align*}
\|\tilde\bdx-\bdx_*\|_2 &~~\le~~
\|(\tilde\bdx_1 +\tilde{V}_2\,\tilde{V}_2^\h\,\bdx_1)-\bdx_1\|_2
+\|\tilde{V}_2\,\tilde{V}_2^\h\,\bdx_2-\bdx_2\|_2 
\end{align*}
while, similar to the proof of Theorem~\ref{t:ls1s} from (\ref{DlA}),
\[\|\tilde{V}_2\,\tilde{V}_2^\h\,\bdx_2-\bdx_2\|_2 ~~\le~~
\frac{\sg_1(A)}{\sg_r(A)}\, 
\frac{1}{1-\frac{\sg_{r+1}+\|\Dl A\|_2}{\sg_r(A)}}\,
\frac{\|\Dl A\|_2}{\|A\|_2}\,2\,\|\tilde{V}_2\,\tilde{V}_2^\h\,\bdx_2\|_2 
\]
leading to (\ref{ics}). ~\Qed

\bibliographystyle{siamplain}

\newpage  \setcounter{page}{1}
\centerline{\Large Online Supplement to} 
\vspace{2mm}
\centerline{\Large ``On the Sensitivity of Singular and Ill-Conditioned }
\vspace{2mm}
\centerline{\Large Linear Systems''}

\vspace{2mm}
\centerline{Zhonggang Zeng\footnote{Department of Mathematics,
Northeastern Illinois University, Chicago, Illinois 60625, USA.
~~email:~{\tt zzeng@neiu.edu}. 
~Research is supported in part by NSF under grant DMS-1620337.a}}

\vspace{2mm}
{\footnotesize
{\bf Abstract.} ~This online supplement provides a software demo of the
package {\sc NAClab} and additional computing examples for calculating
numerical solutions of singular and ill-conditioned linear systems.
}

\renewcommand\thesection{\arabic{section}}

\setcounter{section}{0}

\vspace{2mm}
In this online supplementary material, we briefly introduce the software
package {\sc NAClab} in the context of solving singular linear systems
for the general numerical solution elaborated in the paper 
{\em On the sensitivity of singular and ill-conditioned linear systems}. 
~All the example numbers point to the examples in the paper and all the 
citation numbers point to the references of the paper.

\section*{1. {\sc NAClab} functionality {\tt LinearSolve}}

{\tt NAClab}\footnote{\tt http://homepages.neiu.edu/$\sim$zzeng/naclab.html}
is a software package of Matlab functions for 
numerical algebraic computation \cite{naclab}. 
~We implemented the computation of general numerical solution 
~$\sol_\theta(A,\bdb)$ ~as a functionality {\tt LinearSolve} \cite{solve}
in a simple call with input ~$A$, ~$\bdb$ ~and ~$\theta$:

\vspace{2mm}
\begin{verbatim}
 >> [x0, N, lcnd, res] = LinearSolve(A, b, theta)
\end{verbatim}
\vspace{2mm}

The output {\tt x0}, ~{\tt N}, ~{\tt lcnd} ~and ~{\tt res} carries 
~$A_\theta^\dagger\bdb$, ~$\cK(A_\theta)$ ~spanned by the orthonormal columns 
of ~$N$, ~the sensitivity estimate ~$\|A_\theta\|_2\,\|A_\theta^\dagger\|_2$
~and the residual ~$\max\{\|A\,\bdx_0-\bdb\|_2, ~\|A\,N\|_2\}$ ~respectively.
~The functionality {\tt LinearSolve} follows from the high-rank case of
the template in \S\ref{s:gns}.

Furthermore, {\tt LinearSolve} provide a mechanism to solve a linear system 
in the form of 
\[   L(\bdu_1,\ldots,\bdu_m) ~~=~~ (\bdb_1,\ldots,\bdb_n) ~~~~\mbox{for}~~~~ 
(\bdu_1,\ldots,\bdu_m) 
\]
where ~$L~:~\cU_1\times\cdots\times\cU_m \,\rightarrow\,
\cV_1\times\cdots\times\cV_n$ ~is a linear transformation and 
~$\cU_1,\,\cdots,\,\cU_m,\,\cV_1,\,\cdots,\,\cV_n$ ~are vector spaces of
column vectors, matrices, or polynomials in a call syntax

\vspace{2mm}
\begin{verbatim}
 >> [x0, N, lcnd, res] = LinearSolve({L, domain, parameter}, b, theta)
\end{verbatim}
\vspace{2mm}

\noindent where ~{\tt L} is a Matlab (anonymous) function for evaluating
the linear transformation ~$L$ ~along with ~{\tt domain} ~and 
~{\tt parameter} in cell arrays representing the domain 
~$\cU_1\times\cdots\times\cU_m$ ~and parameters of ~$L$.

\section*{2. Supplement to Example \ref{e:mul}}

Consider a system of polynomial equations ~$\bdf(x,y)\,=\,\bdo$ ~that is
known through empirical data in the perturbed system 
~$\tilde\bdf(x,y)\,=\,\bdo$ ~with
\[ \tilde\bdf(x,y) ~~=~~
\left[ \begin{array}{l}
x^3 + y - 0.7698 \\
x + y^3 - 0.7698 
\end{array} \right]
\]
and the coefficientwise error bound ~$5\times 10^{-7}$.
~A multiple zero 
\[ (x_*,y_*) ~~\approx ~~(\tilde{x},\tilde{y})~~=~~(0.57735,\,0.57735) 
\]
is computed using the {\sc NAClab} polynomial system solver {\tt psolve} and 
the depth-deflation method \cite{DLZ} with an error bound 
~$\|(x_*-\tilde{x},\,y_*-\tilde{y})\|_2\,\le\,\eps\,=\, 9.46\times 10^{-6}$

As briefly elaborated in Example~\ref{e:mul} and in \cite{DLZ}, 
the multiplicity structure can be computed via solving a sequence of 
homogeneous linear systems ~$S_\al(x_*,y_*)\,\bdc\,=\,\bdo$ ~for ~$\al\,=\,
1,2,\ldots$ ~from Macaulay matrices ~$S_\al(\tilde{x},\tilde{y})$ ~serving
as empirical data.
~The {\sc NAClab} functionality {\tt MacaulayMatrix} is built for generating
the Macaulay matrices.
~To construct, say ~$S_2(\tilde{x},\tilde{y})$, ~use the following statements:

\setstretch{0.75}{\scriptsize \noindent
\newline $~~~~${\verb|>> f = {'x^3+y-0.7698','x+y^3-0.7698'};  |}~~\blue{\tt \% cell array of polynomials in character strings}
\newline $~~~~${\verb|>> var = {'x','y'};                      |}~~\blue{\tt \% cell array of variable names in character stings}
\newline $~~~~${\verb|>> z = [.57735,.57735];                  |}~~\blue{\tt \% approximate zero}
\newline $~~~~${\verb|>> M = MacaulayMatrix(z, f, var, 2);     |}~~\blue{\tt \% generate the Macaulay matrix}
\newline $~~~~${\verb|>> single(full(M))                       |}~~\blue{\tt \% display the matrix in single precision}
\newline $~~~~${\verb||}
\newline $~~~~${\verb|ans =|}
\newline $~~~~${\verb||}
\newline $~~~~${\verb|           0   1.0000000   0.9999990           0           0   1.7320499|}
\newline $~~~~${\verb|           0   0.9999990   1.0000000   1.7320499           0           0|}
\newline $~~~~${\verb|           0           0           0   1.0000000   0.9999990           0|}
\newline $~~~~${\verb|           0           0           0   0.9999990   1.0000000           0|}
\newline $~~~~${\verb|           0           0           0           0   1.0000000   0.9999990|}
\newline $~~~~${\verb|           0           0           0           0   0.9999990   1.0000000|}
\newline $~~~~${\verb||}
} 

\setstretch{1.00}{ 
It is a straightforward to verify that the entrywise error on 
~$S_2(\tilde{x},\tilde{y})$ ~is bounded by ~$6\,\eps$ ~on 8 entries.
~As a result, we have an error bound for the error tolerance}
\[  \big\|S_2(x_*,y_*)-S_2(\tilde{x},\tilde{y})\big\|_2
~\le~\big\|S_2(x_*,y_*)-S_2(\tilde{x},\tilde{y})\big\|_{_F} 
~\le~ \sqrt{8}\,6\cdot 9.46\times 10^{-6} ~\approx~ 1.6\times 10^{-4}.
\]
Set the error tolerance slightly larger at
\[ \theta ~~=~~ 2\times 10^{-4} ~~=~~ 0.0002.
\]
Then a one-line call of {\tt LinearSolve} produces the matrix ~$N$ ~whose
columns form an orthonoral basis for the numerical solution 
~$\sol_\theta(S_2(\tilde{x},\tilde{y}),\,\bdo)$ ~in the Grassmannian 
~$\cG_3(\C^6)$.

\setstretch{0.75}{\scriptsize \noindent
\newline $~~~~${\verb|>> [z,N,lcnd,res] = LinearSolve(M,zeros(6,1),2e-4); |}~~\blue{\tt \% solve M*z = 0 within 2e-4}
\newline $~~~~${\verb|>> single(N)                                         |}~~\blue{\tt \% display solution basis in single precision}
\newline $~~~~${\verb||}
\newline $~~~~${\verb|ans =|}
\newline $~~~~${\verb||}
\newline $~~~~${\verb|   1.0000000  -0.0000000   0.0000000|}
\newline $~~~~${\verb|           0  -0.7828174   0.2320854|}
\newline $~~~~${\verb|           0   0.5924006   0.5618970|}
\newline $~~~~${\verb|           0   0.1099370  -0.4584060|}
\newline $~~~~${\verb|           0  -0.1099372   0.4584060|}
\newline $~~~~${\verb|           0   0.1099374  -0.4584059|}
\newline $~~~~${\verb||}
}

\setstretch{1.00}{ 
The multiplicity of ~$(x_*,y_*)$ ~is thus 3, with the dual space
~$\cD_{\bdf,(x_*,y_*)}$ ~accurately represented by the basis}
\begin{align*}
1, ~~&
.78282\,\partial_x +   .59240\,\partial_y +  
.10994\,\mbox{$\frac{1}{2!}$}\,\partial_{x^2}  
-.10994\,\partial_{x\,y}  
+ .10994\,\mbox{$\frac{1}{2!}$}\,\partial_{y^2}\\
& .23209\,\partial_x +
   .56190\,\partial_y 
  -.45840\,\mbox{$\frac{1}{2!}$}\,\partial_{x^2} +
   .45840\,\partial_{x\,y} 
  -.45841\,\mbox{$\frac{1}{2!}$}\,\partial_{y^2}
\end{align*}
so that those differential operators vanish on the entire ideal generated
by the polynomial system at the zero point ~$(x_*,y_*)$.

\section*{3. Supplement to Example~\ref{e:syl}}

The linear equation ~$A(t)\,X+X\,B(t)\,=\,C(t)$ ~can be written as
~$L(X)\,=\,C(t)$ ~where ~$L$ ~is a linear transformation
\[  \begin{array}{ccrcl}
L &~:~& \C^{2\times 2} & ~~\longrightarrow~~& \C^{2\times 2} \\ 
& & X & \longmapsto & A(t)\,X+X\,B(t)
\end{array}
\]
with the domain is ~$\C^{2\times 2}$ ~and parameters
~$A(t)$, ~$B(t)$ ~and ~$C(t)$ ~in (\ref{abct}), 
~The linear system can be solved by constructing the representation matrix and
vectors for ~$L$ ~and ~$C(t)$ ~and solving the resulting matrix-vector 
equation.
~{\tt LinearSolve} in {\sc NAClab} provides an intuitive WYSIWYG approach
for solving the equation directly.
~The matrix-vector representation is generated internally. 

At the hypothetical ~$\tilde{t}\,\approx\,0.6666$ ~with an error bound 0.0001 
of a $4\times 4$ ~system, it clearly safe to say the 2-norm data error bound
is ~$\theta\,=\,10\cdot 0.0001\,=\,10^{-3}$ ~that can be used as the error 
tolerance for the general numerical solution.

\setstretch{0.75}{\scriptsize \noindent
\newline $~~~~${\verb|>> L = @(X,t) [1 -1; 1 -1]*X+X*[-5/3+t 1; -1 -1/3+2*t];            |}~~\blue{\tt \% the linear transformation L}
\newline $~~~~${\verb|>> C = [1 0; 2 -1];                                                |}~~\blue{\tt \% the right-hand side}
\newline $~~~~${\verb|>> domain = {ones(2,2)};                                           |}~~\blue{\tt \% domain of 2x2 matrices}
\newline $~~~~${\verb|>> parameter = {0.6666};                                           |}~~\blue{\tt \% parameter t = 0.6666}
\newline $~~~~${\verb|>> [x0, N, lcnd, res] = LinearSolve({L,domain,parameter}, C, 1e-3);|}~~\blue{\tt \% solve L(X)=C within 1e-3}
\newline $~~~~${\verb|x0 = |}
\newline $~~~~${\verb|    [2x2 double]|}
\newline $~~~~${\verb|N = |}
\newline $~~~~${\verb|    {1x1 cell}    {1x1 cell}|}
\newline $~~~~${\verb|lcnd =|}
\newline $~~~~${\verb|   1.573435327501125|}
\newline $~~~~${\verb|res =|}
\newline $~~~~${\verb|     2.499756923490804e-05|}
\newline $~~~~${\verb||}
}

\setstretch{1.00}{ 
The underlying singular linear system is quite well conditioned with a 
sensitivity measure roughly 1.6 along with a residual ~$\approx\,2.5\times 10^{-5}$, ~implying the general numerical solution carried in {\tt x0} and {\tt N}
are as accurate as the data.}

\setstretch{0.75}{\scriptsize \noindent
\newline $~~~~${\verb|>> x0{1}          |}~~\blue{\tt \% display the truncated SVD solution}
\newline $~~~~${\verb|ans =|}
\newline $~~~~${\verb|   0.249983334213952  -0.250004166457633|}
\newline $~~~~${\verb|  -0.750004165972904  -0.249974998284764|}
\newline $~~~~${\verb||}
\newline $~~~~${\verb|>> N{1}{1}        |}~~\blue{\tt \% display the 1st kernel component}
\newline $~~~~${\verb|ans =|}
\newline $~~~~${\verb|  -0.662148424976858   0.483868243696442|}
\newline $~~~~${\verb|  -0.483822831115126   0.305526519527407|}
\newline $~~~~${\verb||}
\newline $~~~~${\verb|>> N{2}{1}        |}~~\blue{\tt \% display the 2nd kernel component}
\newline $~~~~${\verb|ans =|}
\newline $~~~~${\verb|   0.558171384891092   0.126097939805073|}
\newline $~~~~${\verb|  -0.126073928483796   0.810339052016092|}
\newline $~~~~${\verb| |}
}

\setstretch{1.00}{ 
The result is an accurate approximation to the exact solution (\ref{751}).}

\section*{4. Solving the system (\ref{xu})}

The linear system (\ref{xu}) can be considered as the equation
\[  L(X,U) ~~=~~ (E,\,-F)
\]
where ~$L$ ~is the linear transformation
\[ \begin{array}{ccrcl}
L & : & \C^{3\times 2}\times\C^{1\times 2} & ~~\longrightarrow~~ & 
\C^{3\times 2}\times\C^{1\times 2} \\
& & (X,U) & \longmapsto & (X\,A-B\,X-C\,U, ~D\,X) \end{array}
\]
along with parameters
\begin{align*}
A & ~=~
\mbox{\scriptsize $\left[\begin{array}{cc}1 & 1 \\ 0 & 1 \end{array}\right]$},
~~~
B ~=~ \mbox{\scriptsize $\left[\begin{array}{ccc} 0 & 1 & 0 \\ 0 & 0 & 1 \\ 
2 & -1 & 0 \end{array}\right]$}, ~~~
C ~=~ 
\mbox{\scriptsize $\left[\begin{array}{c} 0 \\ 0 \\ 1 \end{array}\right]$}, ~~~
D ~=~ 
\mbox{\scriptsize $\left[\begin{array}{ccc} 1 & 0 & -1 
\end{array}\right]$},\\
E & ~=~ 
\mbox{\scriptsize $\left[\begin{array}{rr}2 & 1 \\ -1 & 1 \\ 0 & 0 \end{array}
\right]$}~~~
F ~=~
\mbox{\scriptsize $\left[\begin{array}{rr} -1 & 0\end{array}
\right]$}
\end{align*}

In preparation for calling {\tt LinearSolve}, define the linear transformation, 
its domain and the parameter array:

\setstretch{0.75}{\scriptsize \noindent
\newline $~~~~${\verb|>> L = @(X,U,A,B,C,D) {X*A-B*X-C*U,  D*X};                   |}~~\blue{\tt \% the linear transformation function}
\newline $~~~~${\verb|>> A = [1 1; 0 1];  B = [0 1 0; 0 0 1; 2 -1 0]; C = [0; 0; 1]; D = [1 0 -1]; |}~~\blue{\tt \% parameters A,B,C,D}
\newline $~~~~${\verb|>> E = [2 1; -1 1; 0 0];  F = [-1 0];                                        |}~~\blue{\tt \% right side E and F}
\newline $~~~~${\verb|>> domain = {ones(3,2), ones(1,2)};                              |}~~\blue{\tt \% domain of 3x2 and 1x2 matrices}
\newline $~~~~${\verb|>> parameter = {A,B,C,D};                                                  |}~~\blue{\tt \% parameter cell array}
\newline $~~~~${\verb||}
}

\setstretch{1.00}{
The data are exact but floating point arithmetic will introduce entrywise
error arround the unit roundoff ~$\eps\,\approx\,2.2\times 10^{-16}$ ~so that
we can set the error tolerance slightly larger, say ~$\theta\,=\,10^{-10}$.
A simple call to find the general numerical solution within the error
tolerance:

\setstretch{0.75}{\scriptsize \noindent
\newline $~~~~${\verb|>> [Z, N, lcnd, res] = LinearSolve({L,domain,parameter}, {E, -F}, 1e-10) |}~~\blue{\tt \% solve L(X,U) = (E,-F)}
\newline $~~~~${\verb|Z = |}
\newline $~~~~${\verb|    [3x2 double]    [1x2 double]|}
\newline $~~~~${\verb|N = |}
\newline $~~~~${\verb|    {1x2 cell}|}
\newline $~~~~${\verb|lcnd =|}
\newline $~~~~${\verb|  11.987437447750866|}
\newline $~~~~${\verb|res =|}
\newline $~~~~${\verb|     1.332267629550188e-15|}
}

\setstretch{1.00}{
The sensitivity measure approximately 11.99 along with the residual 
~$1.33\times 10^{-15}$ ~indicate that the numerical solution is accurate.}

\setstretch{0.75}{\scriptsize \noindent
\newline $~~~~${\verb|>> Z{1}    |}~~\blue{\tt \% display X component of the 
truncated SVD solution}
\newline $~~~~${\verb|ans =|}
\newline $~~~~${\verb|   1.999999999999999  -0.333333333333333|}
\newline $~~~~${\verb|                   0   0.666666666666666|}
\newline $~~~~${\verb|   0.999999999999999  -0.333333333333333|}
\newline $~~~~${\verb||}
\newline $~~~~${\verb|>> Z{2}    |}~~\blue{\tt \% display U component of the 
truncated SVD solution}
\newline $~~~~${\verb|ans =|}
\newline $~~~~${\verb|  -2.999999999999999   1.999999999999998|}
\newline $~~~~${\verb||}
\newline $~~~~${\verb|>> N{1}{1}    |}~~\blue{\tt \% display the X component of the kernel basis|}
\newline $~~~~${\verb|ans =|}
\newline $~~~~${\verb|                   0  -0.577350269189626|}
\newline $~~~~${\verb|                   0  -0.577350269189626|}
\newline $~~~~${\verb|                   0  -0.577350269189626|}
\newline $~~~~${\verb||}
\newline $~~~~${\verb|>> N{1}{2}    |}~~\blue{\tt \% display the U component of the kernel basis|}
\newline $~~~~${\verb|ans =|}
\newline $~~~~${\verb|     0     0|}
\newline $~~~~${\verb||}
}

\setstretch{1.00}{
Namely, the general numerical solution is an accurate approximation to}
\[ (X,U) ~=~ \left( 
\left[\mbox{\scriptsize $\begin{array}{rr}
2 & -\frac{1}{3} \\ 0 & \frac{2}{3} \\ 1 & -\frac{1}{3} \end{array}$}\right],
~
\left[\mbox{\scriptsize $\begin{array}{cc} -3 & 2 \end{array}$}\right]\right) 
+ t\,\left(
\left[\mbox{\scriptsize $\begin{array}{rr}
0 & -\frac{1}{\sqrt{3}} \\ 0 & -\frac{1}{\sqrt{3}} \\ 0 & -\frac{1}{\sqrt{3}} 
\end{array}$}\right],
~
\left[\mbox{\scriptsize $\begin{array}{cc} 0 & 0 \end{array}$}\right]\right) 
\]

\section*{5. Supplement to Example~\ref{e:uf} and Example~\ref{e:uf1}}

The problem of calculating the B\'ezout coefficients in Example~\ref{e:uf}
and Example~\ref{e:uf1} can be written as
\[   L(u_1,\,u_2,\,u_3) ~~=~~ g  ~~~~\mbox{for}~~~~
(u_1,\,u_2,\,u_3)\,\in\,\mP_3\times\mP_1\times\mP_2
\]
through the linear transformation
\[  \begin{array}{ccrcl}
L & : & \mP_3\times\mP_1\times\mP_2 & ~\longrightarrow~ & \mP_8 \\
& & (u_1,\,u_2,\,u_3) & \longmapsto & u_1\,f_1 + u_2\,f_2 + u_3\,f_3  .
\end{array}
\]
{\sc NAClab} provides an interface for handling polynomials as character 
strings in WYSIWYG manner.
~Polynomial parameters are entered as character strings:

\setstretch{0.75}{\scriptsize \noindent
\newline $~~~~${\verb|>> f1 = '2.5714 + 3.8571*x - 3*x^2 - 6.4286*x^3 - 2.1429*x^4'; |}~~\blue{\tt \% polynomials as character strings}
\newline $~~~~${\verb|>> f2 = '-1.7143 -  1.7143*x + 0.4286*x^2 +  0.4286*x^3 - 3.4286*x^5 -  5.1429*x^6 - 1.7143*x^7'; |}
\newline $~~~${\verb|>> f3 = '0.8571 +  1.2857*x + 2.1429*x^2 +  2.5714*x^3 + 3.4286*x^4 + 3.8571*x^5 + 1.2857*x^6'; |}
\newline $~~~${\verb|>> g =  '4.6667 +  7*x + 2.3333*x^2';                          |}~~\blue{\tt \% the known numerical gcd}
}

\setstretch{1.00}{

\vspace{2mm}
The package {\sc NAClab} provides a library of polynomial operation 
functionalities, such
as ~{\tt pplus(...)} for adding polynomials and ~{\tt ptimes(...)} for 
multiplying polynomials, so that the linear transformation can be defined
as a Matlab (anonymous) function:}

\setstretch{0.75}{\scriptsize \noindent
\newline $~~~~${\verb|>> L = @(u1,u2,u3,f1,f2,f3) ...  |}~~\blue{\tt \% linear transformation (u1,u2,u3) -> u1*f1 + u2*f2 + u3*f3}
\newline $~~~~${\verb|         pplus(ptimes(u1,f1),ptimes(u2,f2),ptimes(u3,f3));|} 
}

\setstretch{1.00}{

\vspace{2mm}
To set the error tolerance, consider the entrywise error bound 
~$0.5\times 10^{-4}$ ~on 55 nonzero entries and}
\[  \|A-\tilde{A}\|_2 ~\le~\|A-\tilde{A}\|_{_F} ~\le~
\sqrt{55}\cdot 0.5\times 10^{-4} ~\approx~ 3.8\times 10^{-4}.
\]
Thus the error tolerance can be set slightly larger at, say ~$5\times 10^{-4}$.
~We can then define the domain and parameter cell arrays and execute 
{\tt LinearSolve} to calculate the B\'ezout coefficients:

\setstretch{0.75}{\scriptsize \noindent
\newline $~~~~${\verb|>> domain = {'1+x+x^2+x^3','1+x','1+x+x^2'};   |}~~\blue{\tt \% domain of polynomials of degrees 3, 1, 2}
\newline $~~~~${\verb|>> parameter = {f1, f2, f3}                    |}~~\blue{\tt \% parameter cell array}
\newline $~~~~${\verb|>> [z0,N,lcnd,res] = LinearSolve({L,domain,parameter}, g, 5e-4) |}~~\blue{\tt \% solve L(u1,u2,u3) = g}
\newline $~~~~${\verb|z0 = |}
\newline $~~~~${\verb|    '0.907108855304999 + 0.333222892924586*x + 0.710289197713311*x^2 + 0.599677838683852*x^3'|}      
\newline $~~~~${\verb|'-0.799463013829436 + 0.0669420537219249*x'      '1.12432524246405 - 0.0664832652437786*x|}
\newline $~~~~${\verb| + 0.0892574807423333*x^2'|}
\newline $~~~~${\verb|N = |}
\newline $~~~~${\verb|    {1x3 cell}    {1x3 cell}|}
\newline $~~~~${\verb|lcnd =|}
\newline $~~~~${\verb|  20.302846223563613|}
\newline $~~~~${\verb|res =|}
\newline $~~~~${\verb|     1.832045500993470e-05|}
}

\setstretch{1.00}{

\vspace{2mm}
The sensitivity is healthy at 20.3 with a residual ~$1.8\times 10^{-5}$ ~so
that the error on the computed general numerical solution is at the same order
of the data error.
~The output {\tt z0} carries the numerical truncated SVD solution 
~$(u_{01},u_{02},u_{03})$ ~as shown above.
~The components ~$(u_{11},u_{12},u_{13})$ ~and ~$(u_{21},u_{22},u_{23})$ ~
of the general numerical solution are in the output ~{\tt N} consists of 
an orthonormal basis of the numerical kernel of the linear transformation
~$L$ ~such as the result shown in Example~\ref{e:uf1}.
~Notice that the output in {\tt z0} and ~{\tt N} carries polynomial in 
WYSIWYG style as character strings. 
~There is no need for a reverse representation and interpretation of
a solution in column vectors.
~Computation of the numerical inverse of the polynomial transformation
matrix shown at the end of Example~\ref{e:uf1} is out the scope of this
paper.
}

\section*{6. An application in solving an integral equation with
an annihilator}

Consider a Volterra integral equation of the first kind in the form of
\begin{equation}\label{vol}
\int_0^s k(s-t)\,x(t)\, dt ~~=~~ g(s), ~~~~~0~\le~s~\le~1
\end{equation}
for finding ~$x(t)$ ~on the interval ~$[0,1]$ ~from the given kernel function 
~$k$ ~and the right-hand side function ~$g$ ~defined on the same interval. 
~The equation is singular if there exists an {\em annihilator} ~$\phi(t)$ 
~such that ~$\int_0^s k(s-t)\,\phi(t)\, dt \,\equiv\, 0$ ~for 
~$0~\le~s~\le~1$.
~As described in \cite[pp. 82-83]{HansInv}, the kernel%
\footnote{There is apparently a typo in \cite[p. 83]{HansInv} about the kernel
(\ref{volk}).}
\begin{equation}\label{volk}  
k(\tau) ~~=~~ \frac{\tau^{-\frac{3}{2}}\,e^{-\frac{1}{4\,\kappa^2\,\tau}}}{
2\,\kappa\,\sqrt{\pi}}
\end{equation}
corresponds to an annihilator ~$\dl(t-1)$, ~the delta function at ~$t=1$.
~For integer ~$n\,>\,0$, ~stepsize ~$h\,=\,\frac{1}{n}$ ~and nodes
~$t_i\,=\,\frac{i}{n}$, ~$i\,=\,0,1,\ldots,n$,
~the equation (\ref{vol}) can be discretized by the linear spline 
approximation
\begin{align*}
   x(t) ~~\approx~~ & \left(1-\frac{t-t_{i-1}}{h} \right) z_{i-1} +
\frac{t-t_{i-1}}{h}  z_i  \\
& ~~~\mbox{for}~~~ t_{i-1}~\le~t~\le~t_i, ~~~~
i~=~ 1, 2, \ldots, n
\end{align*}
and represented by a linear system ~$A\,\bdz\,=\,\bdb$ ~where the variable
~$\bdz\,=\,[z_0,\,\ldots,\,z_n]^\top$ ~in ~$\C^{n+1}$, ~the right-hand side 
vector ~$\bdb\,=\,[b_1,\ldots,b_n]^\top\,\in\,\C^n$ ~and the coefficient 
matrix ~$A\,=\,[a_{ij}]
\,\in\,\C^{n\times (n+1)}$ ~with
\begin{align*}
a_{i1} & ~=~ \int_{t_0}^{t_1} k(t_i-t)\,\left(1-\frac{t-t_0}{h}\right)\,dt \\
a_{ij} & ~=~ 
\int_{t_{j-1}}^{t_j} k(t_i-t)\,\frac{t-t_{j-1}}{h}\,dt +
\int_{t_j}^{t_{j+1}} k(t_i-t)\,\left(1-\frac{t-t_j}{h}\right)\,dt \\
& ~~~~~~~~~~~~~~~~~~~~~~~~~~~~~~~~ j \,=\, 1, 2, \ldots, i+1 \\
a_{i,i+1} & ~=~ \int_{t_{i-1}}^{t_i} k(t_i-t)\,\frac{t-t_{i-1}}{h}\,dt \\
b_i &~=~ g(t_i) \\
&~~~~~\mbox{for}~~~ i = 1, 2, \ldots, n.
\end{align*}
As an experiment with the kernel (\ref{volk}) where ~$\kappa\,=\,4$, ~the
right-hand side function
\[  g(s) ~~=~~ \int_0^s k(s-t)\,dt
\]
of the equation (\ref{vol}) corresponds to a known general solution 
\[ x(t) ~~=~~ 1 + c\,\dl(t-1)
\]
where ~$c$ ~is an arbitrary constant.
~By any standard numerical integration method such as the composite Simpson's 
rule, the matrix ~$A$ ~and the right-hand side vector ~$\bdb$ ~can be 
generated easily. 

Notice that the system ~$A\,\bdz\,=\,\bdb$ ~is underdetermined with the size
~$n\times (n+1)$ ~in addition to being empirical data of a singular equation 
(\ref{vol}).
~We choose not to add an extra equation to square the system for the 
consideration of the inherent singularity in the underlying problem.

For ~$n\,=\,1024$, ~a simple call of {\tt LinearSolve} produces the 
truncated SVD solution ~{\tt z}, ~the matrix ~{\tt K} ~whose columns form an 
orthonormal basis for ~$\cK(A_\theta)$, ~the sensitivity measure ~{\tt lcnd}
~defined as ~$\|A_\theta\|_2\,\|A_\theta^\dagger\|_2$, ~and the
residual ~{\tt res}:

\setstretch{0.75}{\scriptsize \noindent
\newline $~~~~${\verb|>> [z,K,lcnd,res] = LinearSolve(A, b, 1e-6); |}~~\blue{\tt \% solve A*z = b}
\newline $~~~~${\verb|>> lcnd   |}~~\blue{\tt \% the sensitivity}
\newline $~~~~${\verb|lcnd = |}
\newline $~~~~${\verb|     2.469428269074639e+04 |}
\newline $~~~~${\verb| |}
\newline $~~~~${\verb|>> res    |}~~\blue{\tt \% the residual}
\newline $~~~~${\verb|res = |}
\newline $~~~~${\verb|     8.371530784514738e-10 |}
\newline $~~~~${\verb| |}
}

\setstretch{1.00}{ }
The numerical kernel ~$\cK(A_\theta)$ ~is of dimension three. 
~Figure~\ref{f:vol} shows the plot of the truncated SVD solution
~$\bdz$ ~and the numerical kernel vectors ~$\bdu$, ~$\bdv$ ~and ~$\bdw$.

\begin{table}[ht]
\begin{center}
\begin{tabular}{||r|r||r||rrr||} \hline\hline
\multicolumn{1}{||c|}{ }   & 
\multicolumn{1}{c||}{ }    & \multicolumn{1}{c||}{trunc. SVD} & 
\multicolumn{3}{c||}{ basis for ~$\cK(A_\theta)$} \\  
\multicolumn{1}{||c|}{$i$}   & 
\multicolumn{1}{c||}{$t_i$}    & \multicolumn{1}{c||}{solution ~$z_i$} & 
\multicolumn{1}{c}{$u_i$}  &  \multicolumn{1}{c}{$v_i$}  &  
\multicolumn{1}{c||}{$w_i$}  \\
\hline
\tt \scriptsize 0 & \tt \scriptsize         0 &\tt \scriptsize 1.0000025 &\tt \scriptsize  0.0000000 &\tt \scriptsize -0.0000000 &\tt \scriptsize -0.0000000\\
\tt \scriptsize 1 & \tt \scriptsize 0.0009766 &\tt \scriptsize 0.9999977 &\tt \scriptsize -0.0000000 &\tt \scriptsize  0.0000000 &\tt \scriptsize  0.0000000\\
\tt \scriptsize 2 & \tt \scriptsize 0.0019531 &\tt \scriptsize 1.0000017 &\tt \scriptsize -0.0000000 &\tt \scriptsize -0.0000000 &\tt \scriptsize -0.0000000\\
\tt \scriptsize 3 & \tt \scriptsize 0.0029297 &\tt \scriptsize 0.9999994 &\tt \scriptsize  0.0000000 &\tt \scriptsize  0.0000000 &\tt \scriptsize -0.0000000\\
\tt \scriptsize 4 & \tt \scriptsize 0.0039063 &\tt \scriptsize 1.0000004 &\tt \scriptsize -0.0000000 &\tt \scriptsize -0.0000000 &\tt \scriptsize  0.0000000\\
$\cdots$  &$\cdots$~~ &$\cdots$~~ &$\cdots$~~ &$\cdots$~~ & $\cdots$~~ \\
\tt \scriptsize 1007 & \tt \scriptsize 0.9833984 &\tt \scriptsize 1.0000002 &\tt \scriptsize -0.0000001 &\tt \scriptsize -0.0000000 &\tt \scriptsize  0.0000000\\
\tt \scriptsize 1008 & \tt \scriptsize 0.9843750 &\tt \scriptsize 0.9999998 &\tt \scriptsize  0.0000004 &\tt \scriptsize  0.0000001 &\tt \scriptsize -0.0000000\\
\tt \scriptsize 1009 & \tt \scriptsize 0.9853516 &\tt \scriptsize 1.0000010 &\tt \scriptsize -0.0000011 &\tt \scriptsize -0.0000002 &\tt \scriptsize  0.0000000\\
\tt \scriptsize 1010 & \tt \scriptsize 0.9863281 &\tt \scriptsize 0.9999978 &\tt \scriptsize  0.0000032 &\tt \scriptsize  0.0000004 &\tt \scriptsize -0.0000000\\
\tt \scriptsize 1011 & \tt \scriptsize 0.9873047 &\tt \scriptsize 1.0000067 &\tt \scriptsize -0.0000093 &\tt \scriptsize -0.0000013 &\tt \scriptsize  0.0000000\\
\tt \scriptsize 1012 & \tt \scriptsize 0.9882813 &\tt \scriptsize 0.9999813 &\tt \scriptsize  0.0000268 &\tt \scriptsize  0.0000037 &\tt \scriptsize -0.0000001\\
\tt \scriptsize 1013 & \tt \scriptsize 0.9892578 &\tt \scriptsize 1.0000541 &\tt \scriptsize -0.0000770 &\tt \scriptsize -0.0000106 &\tt \scriptsize  0.0000003\\
\tt \scriptsize 1014 & \tt \scriptsize 0.9902344 &\tt \scriptsize 0.9998449 &\tt \scriptsize  0.0002211 &\tt \scriptsize  0.0000303 &\tt \scriptsize -0.0000009\\
\tt \scriptsize 1015 & \tt \scriptsize 0.9912109 &\tt \scriptsize 1.0004460 &\tt \scriptsize -0.0006352 &\tt \scriptsize -0.0000870 &\tt \scriptsize  0.0000025\\
\tt \scriptsize 1016 & \tt \scriptsize 0.9921875 &\tt \scriptsize 0.9987195 &\tt \scriptsize  0.0018246 &\tt \scriptsize  0.0002500 &\tt \scriptsize -0.0000072\\
\tt \scriptsize 1017 & \tt \scriptsize 0.9931641 &\tt \scriptsize 1.0036784 &\tt \scriptsize -0.0052409 &\tt \scriptsize -0.0007182 &\tt \scriptsize  0.0000206\\
\tt \scriptsize 1018 & \tt \scriptsize 0.9941406 &\tt \scriptsize 0.9894347 &\tt \scriptsize  0.0150535 &\tt \scriptsize  0.0020628 &\tt \scriptsize -0.0000591\\
\tt \scriptsize 1019 & \tt \scriptsize 0.9951172 &\tt \scriptsize 1.0303428 &\tt \scriptsize -0.0432346 &\tt \scriptsize -0.0059230 &\tt \scriptsize  0.0001698\\
\tt \scriptsize 1020 & \tt \scriptsize 0.9960938 &\tt \scriptsize 0.9129785 &\tt \scriptsize  0.1240610 &\tt \scriptsize  0.0169516 &\tt \scriptsize -0.0004872\\
\tt \scriptsize 1021 & \tt \scriptsize 0.9970703 &\tt \scriptsize 1.2462977 &\tt \scriptsize -0.3529420 &\tt \scriptsize -0.0470151 &\tt \scriptsize  0.0013861\\
\tt \scriptsize 1022 & \tt \scriptsize 0.9980469 &\tt \scriptsize 0.3926255 &\tt \scriptsize  0.9206328 &\tt \scriptsize  0.0892117 &\tt \scriptsize -0.0036175\\
\tt \scriptsize 1023 & \tt \scriptsize 0.9990234 &\tt \scriptsize 0.0127561 &\tt \scriptsize -0.1016541 &\tt \scriptsize  0.9947379 &\tt \scriptsize  0.0003994\\
\tt \scriptsize 1024 & \tt \scriptsize 1.0000000 &\tt \scriptsize 0.0000008 &\tt \scriptsize  0.0039290 &\tt \scriptsize  0         &\tt \scriptsize  0.9999923 \\ \hline\hline
\end{tabular}
\end{center}
\caption{The general numerical solution 
~$\bdz+\al\,\bdu+\bt\,\bdv+\gamma\,\bdw$ ~that approximates the
exact underlying solution ~$x(t) \,=\,1 + c\cdot\dl(t-1)$ ~for the 
equation (\ref{vol})
}\label{t:vol}
\end{table}

\begin{figure}[ht]
\begin{center}
\epsfig{figure=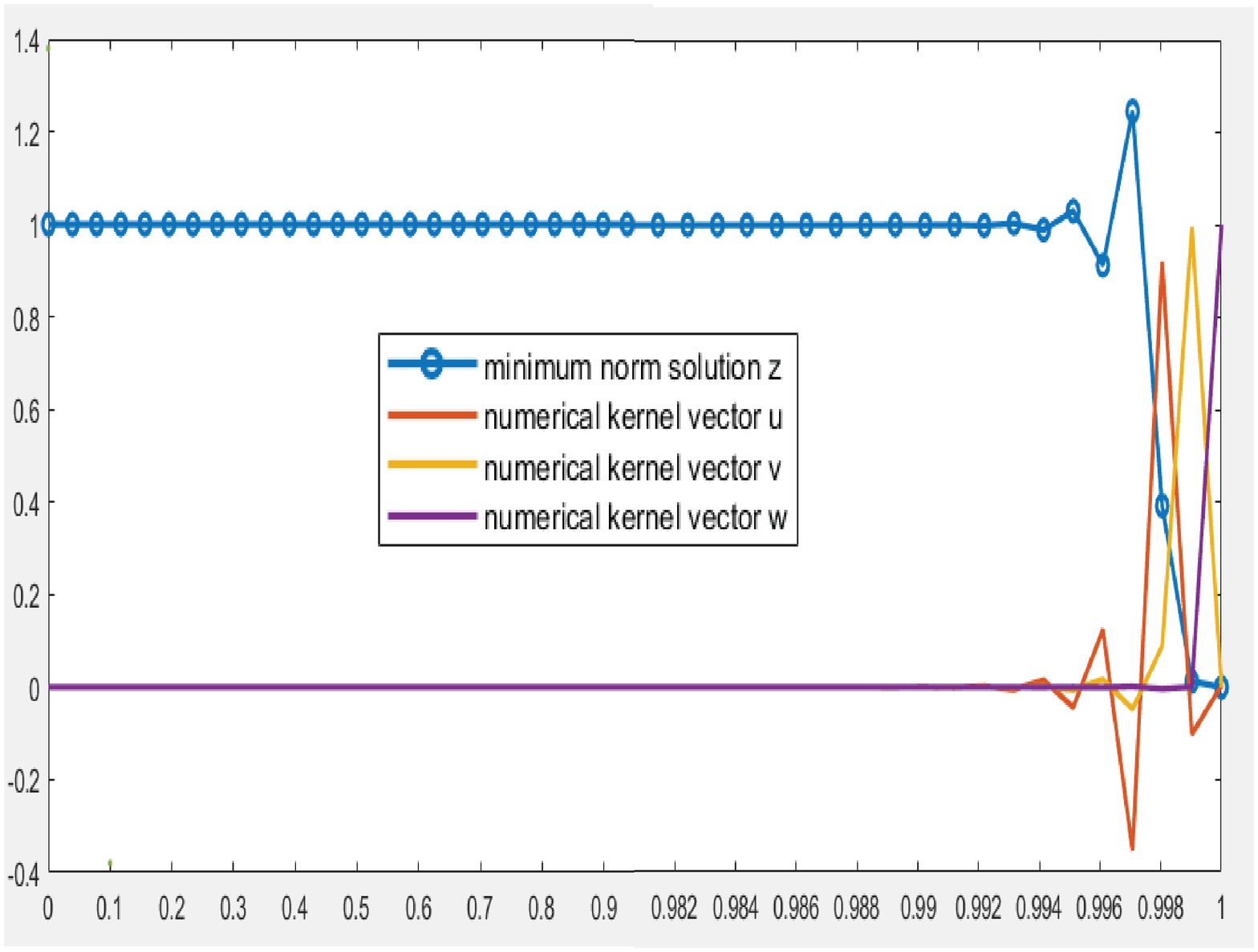,width=5in,height=3in}
\end{center}
\caption{The general numerical solution 
~$\bdz+\al\,\bdu+\bt\,\bdv+\gamma\,\bdw$ ~that approximates the
exact underlying solution ~$x(t) \,=\,1 + c\cdot\dl(t-1)$ ~for the 
equation (\ref{vol})}\label{f:vol}
\end{figure}

Table~\ref{t:vol} shows the actual digits in single precision of the general
numerical solution.
~The condition number
\[ \kappa(A) ~~=~~ \frac{\sg_1(A)}{\sg_n(A)} ~~\approx~~ 3.3\times 10^{25}
\]
is huge whereas the sensitivity of the general numerical solution
~$\bdz+\al\,\bdu+\bt\,\bdv+\gamma\,\bdw$ ~is manageable and much lower
at ~$2.5\times 10^4$.
~Given the residual ~$8.4\times 10^{-10}$, ~we can make a rough error estimate 
of the general numerical solution as
\[  (\|A_\theta\|_2\,\|A_\theta^\dagger\|_2)\,
\left(\|A\,\bdz-\bdb\|_2+\|A\,[\bdu,\,\bdv,\,\bdw]\|_2\right) ~~\approx~~ 
2.1\times 10^{-5}
\]
that can be considered accurate for such an application.

The accuracy estimate can also be justified as follows. 
~The obvious particular solution ~$x_0(t)\,=\,1$ ~of the equation (\ref{vol})
can be approximated by a numerical particular solution:

\setstretch{0.75}{\scriptsize \noindent
\newline $~~~~${\verb|>> y = K\(1-z); |}~~\blue{\tt \% solve K*y+z = 1}
\newline $~~~~${\verb|y = |}
\newline $~~~~${\verb|   0.561958104442570 |}
\newline $~~~~${\verb|   1.049493278421724 |}
\newline $~~~~${\verb|   0.997798939791854 |}
\newline $~~~~${\verb| |}
\newline $~~~~${\verb|>> norm(Z+K*u-1,1)*h  |}~~\blue{\tt \% error of the numerical particular solutoin}
\newline $~~~~${\verb|ans = |}
\newline $~~~~${\verb|     1.090344305094233e-07 |}
\newline $~~~~${\verb| |}
}

\setstretch{1.00}{ }
Namely, the particular solution ~$x_0(t)\,=\,1$ ~can be accurately approximated
by a numerical particular solution with error bound ~$1.1\times 10^{-7}$.
~On the other hand, the solution to the homogeneous equation corresponding to
(\ref{vol}) is
\[ \dl(t-1) ~~=~~ \lim_{\eps\rightarrow 0+}\,\dl_\eps (t-1)
~~~~\mbox{where}~~~~ \dl_\eps (t-1) ~=~
\frac{e^{-\left(\frac{t-1}{\eps}\right)^2}}{\eps\,\sqrt{\pi}}
\]

The following Matlab statement sequence shows that a numerical kernel vector
approximates ~$\dl_{0.00001}(t-1)$ ~with an error measure ~$7.76\times 10^{-5}$,
~as the error estimate suggests.

\setstretch{0.75}{\scriptsize \noindent
\newline $~~~~${\verb|>> y = K\(1-z); |}~~\blue{\tt \% solve K*y+z = 1}
\newline $~~~~${\verb|y = |}
\newline $~~~~${\verb|>> epsilon = 1e-5; |}~~\blue{\tt \% a tiny epsilon}
\newline $~~~~${\verb|>> h = 1/n; |}~~\blue{\tt \% stepsize}
\newline $~~~~${\verb|>> t = 0:h:1; |}~~\blue{\tt \% the nodes in t}
\newline $~~~~${\verb|>> v = K\(1./(epsilon*sqrt(pi)*exp(((t-1)/epsilon).^2))) |}~~\blue{\tt \% solve K*v = delta\_epsilon}
\newline $~~~~${\verb|v = |}
\newline $~~~~${\verb|   1.0e+04 * |}
\newline $~~~~${\verb|   0.022167287609456 |}
\newline $~~~~${\verb|   0.000000000000000 |}
\newline $~~~~${\verb|   5.641852287123326 |}
\newline $~~~~${\verb| |}
\newline $~~~~${\verb|>> norm(K*v -1./(epsilon*sqrt(pi)*exp(((t-1)/epsilon).^2)),1)*h |}
\newline $~~~~${\verb|ans = |}
\newline $~~~~${\verb|     7.764167138103457e-05 |}}

\setstretch{1.00}{ }

\vspace{2mm}
Although the underlying system is underdetermined in addition to being 
singular, the general numerical solution 
~$\bdz+\al\,\bdu+\bt\,\bdv+\gamma\,\bdw$ ~accurately reveals that the equation%
~(\ref{vol}) with the specific kernel (\ref{volk}) can be accurately solved
in an interval ~$[0,1-\eps)$ ~for a small ~$\eps\,>\,0$ ~since the annihilator 
represented by ~$\al\,\bdu+\bt\,\bdv+\gamma\,\bdw$ ~is identically zero in 
that interval.
~The singularity of the underlying equation compounded by the representation 
linear system being underdetermined is not detrimental at all if we compute
the general numerical solution and, in particular, consider the numerical
kernel as an integral part of the solution.

\end{document}